\documentclass[a4paper,12pt,reqno]{amsart}
\usepackage{amssymb, amscd,enumerate,enumitem,xcolor, mathtools, graphics, graphpap}
\usepackage[british]{babel}
\usepackage[mathscr]{euscript}
\usepackage[makeroom]{cancel}
\usepackage[normalem]{ulem}
\usepackage[colorlinks,
hypertexnames=false]{hyperref} 

\setlist[itemize]{itemsep=2pt plus 4pt minus 2pt, topsep=1pt, leftmargin=18pt}

\makeatletter
\def\subsection{\@startsection{subsection}{3}%
  \z@{.5\linespacing\@plus.7\linespacing}{.1\linespacing}%
  {\normalfont\bfseries}}
\makeatother

\font\smallsmc = cmcsc9
\font\smalltt = cmtt8
\font\smallit = cmti8

\numberwithin{equation}{section}

\theoremstyle{plain}
\newtheorem{theo}{Theorem}[section]
\newtheorem{lem}[theo]{Lemma}
\newtheorem{prop}[theo]{Proposition}

\newtheorem{cor}[theo]{Corollary}

\theoremstyle{definition}
\newtheorem{rem}[theo]{Remark}

\newtheorem{definition}[theo]{Definition}

\newenvironment{pf}{\noindent{\it Proof. }}{$\hfill\square$\par\medskip}

\theoremstyle{plain}

\theoremstyle{definition}

\newcommand{\beq}{\begin{equation}}
\newcommand{\eeq}{\end{equation}}
\newcommand{\beqn}{\begin{equation*}}
\newcommand{\eeqn}{\end{equation*}}
\renewcommand{\a}{\alpha}
\renewcommand{\b}{\beta}

\newcommand{\e}{\epsilon}

\newcommand{\f}{\varphi}
\newcommand{\g}{\gamma}
\newcommand{\h}{\eta}

\renewcommand{\k}{\kappa}
\renewcommand{\l}{\lambda}
\renewcommand{\o}{\omega}
\newcommand{\q}{\vartheta}
\renewcommand{\r}{\rho}
\newcommand{\s}{\sigma}

\renewcommand{\u}{\upsilon}
\newcommand{\x}{\xi}
\newcommand{\z}{\zeta}

\newcommand{\G}{\Gamma}

\renewcommand{\L}{\Lambda}
\renewcommand{\O}{\Omega}
\newcommand{\Q}{\Theta}

\newcommand{\bA}{\mathbb{A}}

\newcommand{\bC}{\mathbb{C}}
\newcommand{\bR}{\mathbb{R}}
\newcommand{\bZ}{\mathbb{Z}}

\newcommand{\bH}{\mathbb{H}}

\newcommand{\bI}{\mathbb{I}}
\newcommand{\bJ}{\mathbb{J}}
\newcommand{\bM}{\mathbb{M}}

\newcommand{\gc}{\mathfrak{c}}

\renewcommand{\gg}{\mathfrak{g}}

\newcommand{\gl}{\mathfrak{l}}

\newcommand{\gu}{\mathfrak{u}}

\newcommand{\gX}{\mathfrak{X}}

\newcommand{\su}{\mathfrak{su}}

\newcommand\GL{\mathrm{GL}}
\newcommand\SL{\mathrm{SL}}

\newcommand\SU{\mathrm{SU}}
\newcommand\U{\mathrm{U}}

\newcommand\Sp{\mathrm{Sp}}
\renewcommand\sp{\mathfrak{sp}}
\renewcommand\sl{\mathfrak{sl}}
\newcommand\ggl{\mathfrak{gl}}

\newcommand{\cA}{\mathscr{A}}
\newcommand{\cC}{\mathcal{C}}
\newcommand{\cD}{\mathscr{D}}

\newcommand{\cF}{\mathscr{F}}

\newcommand{\cH}{\mathscr{H}}

\newcommand{\cM}{\mathscr{M}}

\newcommand{\cS}{\mathscr{S}}

\newcommand{\cU}{\mathscr{U}}
\newcommand{\cV}{\mathscr{V}}
\newcommand{\cW}{\mathscr{W}}

\newcommand{\cZ}{\mathscr{Z}}

\newcommand{\Jst}{J_o{}} 

\newcommand{\p}{\partial}

\renewcommand{\square}{\kern1pt\vbox
{\hrule height 0.6pt\hbox{\vrule width 0.6pt\hskip 3pt
\vbox{\vskip 6pt}\hskip 3pt\vrule width 0.6pt}\hrule height0.6pt}\kern1pt}

\DeclareMathOperator\Tr{Tr\;}
\DeclareMathOperator\End{End}

\DeclareMathOperator\Ad{Ad}

\DeclareMathOperator\Id{Id}

\DeclareMathOperator{\Span}{span}

\DeclareMathOperator{\Lie}{Lie}
\DeclareMathOperator{\qk}{qk}

\renewcommand\={:=}

\newcommand{\wt}{\widetilde}
\newcommand{\wh}{\widehat}

\newcommand{\bt}{\begin{theo}\ \ }
\newcommand{\et}{\end{theo}}
\newcommand{\bp}{\begin{prop}\ \ }
\newcommand{\ep}{\end{prop}}
\newcommand{\bc}{\begin{cor}\ \ }
\newcommand{\ec}{\end{cor}}
\newcommand{\bl}{\begin{lem}\ \ }
\newcommand{\el}{\end{lem}}
\newcommand{\bd}{\begin{definition}}
\newcommand{\ed}{\end{definition}}
\newcommand{\n}{\nabla}

\newcommand{\be}{\begin{equation}}
\newcommand{\ee}{\end{equation}}

\def\<#1,#2>{\langle\,#1,\,#2\,\rangle}
\newcommand{\arr}{\begin{array}{rlll}}
\newcommand{\ea}{\end{array}}
\newcommand{\bea}{\begin{eqnarray}}
\newcommand{\eea}{\end{eqnarray}}
\newcommand{\bean}{\begin{eqnarray*}}
\newcommand{\eean}{\end{eqnarray*}}

\newcommand{\ZZ}{{\cZ(N)}}
\newcommand{\ZZM}{{\cZ(M)}}
\newcommand{\HM}{{\cH(M)}}
\newcommand{\HN}{{\cH(N)}}
\newcommand{\HMqk}{{\HM^{\operatorname{qk}}}}
\newcommand{\Hqk}{{\cH^{\qk}}}
\newcommand{\HCMqk}{{\HCM^{\operatorname{qk}}}}
\newcommand{\HCM}{{\cH^\bC(M)}}

\newcommand{\Iden}{\operatorname{I}}

\newcommand{\sH}{\mathbf H}
\newcommand{\sE}{\mathbf  E}

\def\sideremark#1{\ifvmode\leavevmode\fi\vadjust{
\vbox to0pt{\hbox to 0pt{\hskip\hsize\hskip1em
\vbox{\hsize3cm\tiny\raggedright\pretolerance10000
\noindent #1\hfill}\hss}\vbox to8pt{\vfil}\vss}}}
\begin{document}

\title[Instantons on  quaternionic K\"ahler manifolds]
{Hyperk\"ahler cones and instantons on \\ quaternionic K\"ahler manifolds}
\author[C. Devchand]{Chandrashekar Devchand}
\author[M. Pontecorvo]{Massimiliano Pontecorvo}
\author[A. Spiro]{Andrea  Spiro}

\begin{abstract}  
We present a novel approach to the study of Yang-Mills instantons on quaternionic K\"ahler manifolds, based on an extension  of the   harmonic space  method of  constructing instantons on  hyperk\"ahler manifolds. Our results  establish a bijection between local equivalence classes of  instantons on quaternionic K\"ahler manifolds $M$ and equivalence classes of  certain holomorphic maps  on an appropriate $\SL_2(\bC)$-bundle over  the Swann bundle of $M$.  
 \end{abstract}
 
 \thanks{This research was partially supported by  
{\it Ministero dell'Istruzione, Universit\`a e Ri\-cer\-ca} in the framework of the project 
``Real and Complex Manifolds: Geometry, Topo\-logy and  Harmonic Analysis''  
and by GNSAGA of INdAM}


\vskip 4truemm
\maketitle

\section{Introduction}
Let $(M, g)$ be a $4n$-dimensional  quaternionic pseudo-K\"ahler (qk) manifold, by which we mean a pseudo-Riemannian manifold with  holonomy in  $\Sp_1{\cdot} \Sp_{p,q}$, $p + q = n$.   Any such  manifold is automatically Einstein   and  its constant scalar curvature is zero  if and only if it is the quotient of   a   pseudo-hyperk\"ahler (hk) manifold by some discrete group of isometries.  In this paper we are  interested in studying instantons on $M$,    that is to say in vector bundles equipped  with connections,  whose curvature tensor $F$ is  pointwise invariant under the  quaternionic structure of  the tangent space.   They are gauge fields which minimise the Yang-Mills functional 
$YM = ||F|| := \int_M \Tr F\wedge * F$, simply because they saturate a topological bound on $YM$,  a  property with  a  number of crucial consequences both in  physics and mathematics. For an  overview of the importance of  instantons in mathematics,  the reader might take a look at   the beautiful lecture notes  \cite{Uh}.  
\par
 In a previous paper \cite{DPS},  we  have considered  instantons on   hk  manifolds  and  presented a differential geometric formulation of the harmonic space approach to them,   a method that was originally developed in theoretical physics for  studying supersymmetric field theories (see e.g. \cite{gios, gios_book}). 
It is    based on the  notion of {\it harmonic space}  $\HN$  for an hk manifold $(N, g)$, which  is the trivial bundle  
$\HN  = N \times  \SL_2(\bC)$ over $N$ equipped with  a non-product complex structure that can be described as  follows.  Each fibre $\{x\} \times \SL_2(\bC)$ of $\HN$ naturally projects onto the quotient 
$\{x\} \times \SL_2(\bC)/B \simeq \bC P^1$  by the  Borel subgroup of $\SL_2(\bC)$ consisting of  the upper triangular matrices.  Thus  $\HN$ can be considered as  a $B$-bundle   over the twistor space $ Z(N) = N \times \bC P^1$ 
of the hk manifold $N$.  The complex structure  of  $\HN$ is precisely  the  unique complex structure which makes 
$\HN$  a  holomorphic bundle over $Z(N)$.  \par
In   \cite{DPS}  we provide full and complete proofs of the following  correspondence,  upon which  the harmonic space method is founded.  Starting from Yang-Mills data $(E,D)$ on an hk manifold  $(N, h)$, where $E$ is a vector bundle associated to a principal $G$-bundle with  connection $D$,  we  show there that  $(E, D)$   is an instanton  if  and only if the corresponding  pull-back $(E', D')$ over $\HN$ admits a  distinguished class of local trivialisations (= gauges) of $E'$ called {\it analytic gauges},  defined by the vanishing of certain components of the gauge potential  $A'$ of the connection $D'$. 
In these gauges,  all the remaining components of $A'$  can be  uniquely reconstructed  from  just one of them, the  so-called {\it prepotential},  a  holomorphic function taking values in the complexification $\gg^\bC$ of the   Lie algebra   of the structure group $G$.  Further,  {\it any} holomorphic map  satisfying a certain simple first-order equation on an 
$\SL_2(\bC)$-invariant open domain $\cW \subset \HN$ is a prepotential for a (locally defined) instanton on $N$.  Combining these  results  establishes a bijective correspondence between a special class of normalised prepotentials on $\cW$ and  the moduli space of locally defined instantons on  $N$.

In the present paper, we establish a similar bijection  for  instantons on  strict  qk manifolds, i.e. those having non-zero scalar curvature.
The notion of  the {\it Swann bundle} $\pi\colon \cS(M) \to M$ of a (strict) qk manifold $(M, g)$ \cite{Sw} is crucial for our discussion.
This is an $\bH^*/\bZ_2$-bundle over $M$, canonically  determined by  the qk structure and naturally equipped  with a distinguished pseudo-Riemannian
metric $h$ of conical type  and a triple $(J_\a)_{\a = 1, 2, 3}$ of integrable complex structures.
It is thus  an hk  manifold,  the  so-called  {\it hyperk\"ahler cone} of $M$  \cite{PPS} (see also \cite{BCDGVV, Ha,Hi}).
Following \cite{DPS}, we may define its harmonic space $\cH(\cS(M))$ as the trivial  $\SL_2(\bC)$-bundle over $\cS(M)$.
This  naturally fibres over $M$, with   fibres  $\SL_2(\bC) \ltimes \bH^*/\bZ_2$.  We call this principal bundle $\cH(\cS(M)) \to M$ the  {\it harmonic space of the qk manifold  $(M, g)$}.  
Now, a  Yang-Mills field $(E, D)$  on  $(M, g)$ may successively be lifted to the  hk manifold $(\cS(M), h, J_\a)$ and to its harmonic space $\cH(\cS(M))$. Then, applying our previous results on instantons on hk manifolds  \cite{DPS} and imposing an $\bH^*/\bZ_2$-symmetry, we obtain the desired bijection between equivalence classes of certain holomorphic maps  and equivalence classes of locally defined  instantons on a strict qk manifold  $M$. 
\par
This paper is structured as follows. After some preliminaries (\S 2), we provide (in \S 3) our presentation of the Swann bundle and  discuss some of its properties which we shall need. 
In \S 4 and \S 5, we prove the above-described  results on the correspondence between prepotentials and equivalence classes of instantons.
In the final section (\S 6) we use the classical results of Narasimhan and Ramanan \cite{NR} on gauge fields with
compact structure groups to express our solutions in a form reminiscent of the ADHM construction, showing that any 
local equivalence class of instantons  is uniquely associated with a corresponding finite set of  (locally defined) smooth
matrix-valued maps $\l\colon \cU \subset M \to  \text{Mat}_{k \times m}(\bC)$. 
This yields a more compact representation of the instanton prepotentials, which we expect will facilitate the further investigation of the moduli space of  instantons on  compact qk manifolds. 
We conclude with  a discussion in our framework of the explicit  examples of instantons on $\bH P^n$ considered   by  Mamone Capria and Salamon  \cite{MS}. 
\par
\noindent
{\bf Acknowledgement.}
CD thanks Hermann Nicolai and the Max-Planck-Institut f\"ur Gravitationsphysik 
(Albert-Einstein-Institut) for providing an excellent research environment.

\section{Preliminaries on quaternionic K\"ahler  manifolds\\ and  instantons}
In this section, in order to fix our notation and definitions,
we review the fundamental properties of quaternionic K\"ahler  manifolds  and of gauge fields,
which we shall need in this paper.
\subsection{Quaternionic K\"ahler manifolds}
\label{notation}
Let $W$ be   a $4n$-dimensional real vector space.
A  {\it hypercomplex structure\/}  on $W$ is a triple $(I_1, I_2, I_3)$ of endomorphisms
satisfying the multiplicative relations of the imaginary quaternions,  
$I^2_\a = - \Id_W$ and $I_\a I_\b =  I_\g$
for all cyclic permutations  $(\a,\b,\g)$  of $(1,2,3)$.
This notion has the following counterpart in the category of manifolds:
a  {\it hypercomplex structure} on a  $4n$-dimensional  real manifold $M$
is a triple $(J_1, J_2, J_3)$ of 
 integrable complex structures on $M$, with the property   that each triple $(I_\a \= J_\a|_x)$, $x \in M$, 
is a hypercomplex structure on  $T_xM$.\par
These definitions admit the following  generalisations.  Let $Q_W=$ $\Span_\bR(I_1$, $ I_2, I_3) \subset \End(W)$ be
the $3$-dimensional space  of endomorphisms spanned by the three elements $I_\a$ of a hypercomplex structure on $W$.
A {\it  quaternionic pseudo-K\"ahler structure} on $W$ is an  inner product $\langle\cdot, \cdot \rangle$ on $W$  which is
{\it hermitian with respect to   $Q_W$}, meaning  that every  $J \in Q_W$ is skew-symmetric  with respect to  $\langle\cdot, \cdot \rangle$.
The corresponding definition for manifolds is as follows.
\begin{definition}\label{defqk}
A $4n$-dimensional pseudo-Riemannian manifold $(M, g)$ 
of signature  $(4p, 4q)$, $p +q = n$,   is  a {\it  quaternionic pseudo-K\"ahler}  ({\it qk}) {\it manifold}
if it  admits a subbundle $Q \subset \End(TM)$ of quaternionic structures on the tangent spaces satisfying  the following two conditions: 
\begin{itemize}
\item[i)]   each inner product $g_x$ is Hermitian  with respect to the  corresponding quaternionic structure $Q_x$; 
\item[ii)]  the  parallel transport of the Levi-Civita connection $\n$ of $g$  preserves $Q$.
\end{itemize}
Further, $(M, g)$ is a {\it pseudo-hyperk\"ahler\/ (hk) manifold} if there exist three global $\nabla$-parallel sections 
$J_1, J_2, J_3$ of $Q$ (which are thus integrable complex structures) giving a hypercomplex structure on each tangent space of  $M$.
\end{definition}
It is well known that a qk manifold is Einstein and that  its scalar curvature is  zero if and only if it is locally isometric to an hk manifold.
A qk manifold that is not locally  hk, i.e. not Ricci flat,  is called a {\it strict} qk manifold.
\par
\smallskip
Let $(M, g, Q)$
 be a $4n$-dimensional qk manifold.  For each $x \in M$,  an {\it adapted  frame} of $T_x M$ is an  orthonormal frame $(e_i) \subset T_xM$, 
   for which there is  a hypercomplex structure  $(J_1, J_2, J_3)  \subset Q_x$  such that   
\beq \label{adapted} e_{(4 k + 1) +\a} = J_\a e_{4k +1}\ \ ,\quad \text{for}\ k = 0,\ldots, n-1\ \text{and}\ \  \a = 1,2,3\ .\eeq
 The adapted  frames  form a principal bundle  $\pi\colon  O_g(M, Q) \to M$ with  structure group    $\Sp_1{\cdot} \Sp_{p,q}$. It  is   preserved by the Levi-Civita connection and  is  a holonomic reduction  of  the orthonormal frame bundle $\pi'\colon O_g(M) \to M$ of $(M, g)$.\par
 \smallskip
Each adapted frame yields a  convenient identification  of the complexified  tangent space $T^\bC_x M$ with a tensor product of the form $\sH_x \otimes \sE_x$, with  $\sH_x \simeq \bC^2$ and $\sE_x \simeq  \bC^{2n}$, first considered  by Salamon  in  \cite{Sa}.
To see this consider first, as  the standard  representation of $\Sp_1{\cdot}\Sp_{p,q}$ on $\bH^n = \bC^{2n}$,   the one  for   which
the element  $ (\smallmatrix i & 0\\ 0 & -i  \endsmallmatrix )\in \sp_1$ acts on $(\bH^n)^\bC \simeq \bC^2 \otimes \bC^{2n}$ as the  left multiplication by  the imaginary unit $\mathbf i \in \bH$.  Now,  if we use  an   adapted frame $\u = (e_i) \in O_g(M, J_\a)|_x$ to make the identification  $T_x M \simeq \bR^{4n} = \bH^n$, the standard left action of  $\Sp_1 {\cdot} \Sp_{p,q}$ on $\bH^n$ corresponds to  a left action of  $\Sp_1 {\cdot} \Sp_{p,q}$  on $T_x M$. By  the above  assumption, this action is such that   $ (\smallmatrix i & 0\\ 0 & -i  \endsmallmatrix)$ 
acts  on $T_x M$   precisely as the action of  complex structure   $I = \u_*(\mathbf i)$. 
This   $\Sp_1 {\cdot} \Sp_{p,q}$-action  extends by $\bC$-linearity to  the
standard  action of  $\SL_2(\bC) {\cdot} \Sp_n(\bC)$ on
  $(\bH^n)^\bC \simeq \bC^2 \otimes \bC^{2n}$.  Consequently,  the  $\bC$-linear extension $\u\colon (\bH^n)^\bC \to T^\bC_x M$ of the considered adapted frame yields
 an  $\SL_2(\bC) {\cdot} \Sp_n(\bC)$-equivariant isomorphism $T^\bC_x M \simeq \sH_x \otimes \sE_x$  
 with   $\sH_x \simeq \bC^2$ and $\sE_x \simeq \bC^{2n}$.
Any two  such isomorphisms are necessarily related  to each other by the  action  of  an element of  $ \Sp_1{\cdot}\Sp_{p,q} \subset \SL_2(\bC){\cdot} \Sp_n(\bC)$,
which clearly preserves the tensor product structure $\sH_x \otimes \sE_x$.
Thus, modulo the action of  this group, the identification of   $T^\bC_x M$ with $\sH_x \otimes \sE_x$  is  independent of the chosen adapted frame $\u = (e_i)$.  \par
\subsection{Gauge  fields and potentials} \label{section2.2}
Let   $p\colon P \to M$  be a  principal $G$-bundle and  $p^E \colon E = P {\times_{G, \r}} V {\to} M$
 an associated vector bundle with fibre $V \simeq \bR^N$, determined by a faithful linear representation $\r\colon G \to \GL(V)$.\par
 \smallskip
 Following physics terminology,   a  {\it gauge} is a local trivialisation   
 $\f\colon P|_{\cU} \to \cU \times G$ of  $P$ over  an open set $\cU$. We call $\cU$ the   {\it domain} of the   gauge. 
Given  two gauges   $\f, \f'$ on  overlapping domains $\cU$, $\cU'$,   the associated {\it transition function} is a map of the form $(\f' \circ \f^{-1})(x, h) = (x, g_x{\cdot} h)$
for a smooth  map $g\colon \cU \cap \cU' \to G$.
The {\it gauge transformation} corresponding to this transition function is the smooth family  of $G$-automorphisms   $h \mapsto g_x{\cdot} h$, $x \in \cU \cap \cU'$.   \par
\smallskip
Given a  connection $1$-form $\o$ on $P$,  there is  a uniquely associated covariant derivative  $D$ for the sections of $E = P \times_{G, \r} V$,  locally defined as follows.  Consider a gauge  $P|_{\cU}\overset \sim \to \cU \times G$ and use it to identify  each tangent space 
$T_x P$, $x \in \cU$,  with a direct sum $ T_{(x, g)}( \cU \times G) \simeq T_x\cU + T_gG$. Then,    the connection $1$-form $\o$ is pointwise identifiable  with a  sum of  the form  $\o|_{(x, g)} = - A_{(x,g)} +  \o^G_g$, 
where  $\o^G$ is the Maurer-Cartan form of   $G$ and each  map   $A_{(\cdot, g)}\colon \cU \to T^* \cU \otimes \gg$
is a   $\gg$-valued $1$-form, changing $G$-equivariantly with respect to $g$. 
The $1$-form  $A \= A_{(\cdot, e)}\colon \cU \to T^*M \otimes \gg$  is   called the
 {\it potential of $\o$} in the considered gauge.  
 If  $P|_{\cU'} \simeq \cU' \times G$ is a different  gauge on a domain  $\cU'$  that  overlaps  with $\cU$ and  $h \mapsto g_x{\cdot} h$ is
 the  corresponding gauge transformation, then  the    potential of $\o$ in this new gauge  is given by
 \beq \label{changepot} A' =  \Ad_{g^{-1}} A  + g^{-1}   dg\qquad \text{on}\ \cU \cap \cU'\ . \eeq
Now, if   $\s\colon \cU \to E|_{\cU}  \simeq \cU \times V $ is a (local) section  of   $E = P\times_{G, \r} V$ of the form
$\s(x) = (x, s^i(x))$ for a smooth map $(s^i)\colon \cU \to V = \bR^N$, then the covariant derivative of $\s$  along a vector field $X \in \gX(\cU)$  is the section $D_X \s$  of $ E|_{\cU}  \simeq \cU \times V $ defined by
\beq \label{covder} D_X\s(x) =  \left(x, (X {\cdot} \s^i)|_x + A_j^i(X) \s^j|_x\right)\ ,\eeq
 where $A^i_j \= \r_*{\circ}A$,  with   $\r_*\colon \gg \to \ggl(V)$   being
the Lie algebra representation determined by  the linear   representation $\r\colon G \to \GL(V)$.  
 \par
\medskip
Given a Lie group $G$, we shall call a pair $(E, D)$  consisting of 
\begin{itemize}[itemsep=1pt, topsep=0pt , leftmargin=18pt]
\item[1)] a vector bundle  $E$,  associated  with   a principal $G$-bundle $p\colon P \to M$ and
\item[2)] a covariant derivative $D$ on $E$, induced by a connection $\o$ on $P$
\end{itemize}
a {\it gauge field with structure group $G$}
and we  shall   say that   {\it  $(E, D)$ is  associated with   $(P, \o)$}.
Consider a complex Lie group $G$  with  compact real form  $G^o \subset G$.
If  $p:P \to M$ and $\o$  admit  reductions to a $G^o$-subbundle $P^o \subset P$   and  a connection $\o^o$ on $P^o$,
we  say that {\it  the pair $(P, \o)$ is reducible to $(P^o, \o^o)$} and that  the gauge field {\it $(E, D)$ is
the {\rm complexification} of a gauge field with compact structure group $G^o$}.\par
\smallskip
 In a gauge $\f\colon P|_{\cU} \to \cU \times G$ and its  associated  gauge $\wh \f$  for  $E = P \times_{G, \r} V$,
$$ \wh \f\colon E|_{\cU} \to \cU \times V\ ,\qquad \wh \f\left([\f^{-1}(x, e), v)]\right) \= (x, v)\ ,$$  the  {\it curvature tensors}
of  $\o$  and   $D$ are respectively the ($\gg$- and $\ggl(V)$-valued)  $2$-forms $\cF^\f$ and $F_x$ on $\cU$,  defined   by
\beq \label{2.4}
\cF^\f_x(v,w) \= 2 \O_{(x,e)}(v, w) \ ,\quad F_x(v, w) \=   \r_* \circ \cF^\f_x(v,w)\ , \quad\text{for }\  v, w \in T_x \cU\ ,\eeq
where $\O$ is the curvature $2$-form of the connection $1$-form $\o$.
Note that  $\cF^\f$ is expressed in terms of  the potential $A$   by 
\beq \label{Fcurvature} \cF^\f(X, Y) = X {\cdot}(A(Y)) - Y {\cdot}(A(X)) + [A(X), A(Y)] - A([X, Y])\ .\eeq
The curvature $F$ of $D$ is independent of the considered gauge and it 
satisfies the identity 
\beq F(X, Y)s  = [D_X, D_Y]s - D_{[X,Y]} s\qquad \text{for sections}\ s \in \G(E)\ ,\eeq
a property that is often used as an alternative  definition for $F$.\par
\smallskip
The above notions  have analogues in  the 
category  of holomorphic bundles.  If $(M, J)$ is a 
complex manifold,  $G$ a complex Lie group and $V$  a complex vector space, then a principal $G$-bundle 
$p\colon P \to M$ (resp. a complex vector bundle $p\colon E \to M$ with fibre $V$) is called {\it holomorphic} if it is equipped with a complex structure $\wh J$, such that the right action of $G$ on $P$ (resp. the vector bundle structure on $E$) is $\wh J$-holomorphic  
and the projection $p$ is a $(\wh J, J)$-holomorphic mapping.  A gauge  is  called  {\it holomorphic} 
if it is a local holomorphic map from $(P, \wh J)$ (resp.  $(E, \wh J)$) to the cartesian product $M \times G$ (resp. $M \times V$), equipped with  the product complex structure.\par
A connection form $\o$ on a holomorphic $G$-bundle $(P, \wh J)$  is called {\it $\wh J$-invariant} if the corresponding horizontal spaces in $TP$ are invariant under the complex structure $\wh J$. This is the case if and only if  in any holomorphic gauge $\f\colon P|_\cU\to \cU \times G$, the potential $A\colon \cU \to  T^* \cU \otimes  \gg$ takes values  in $ T^{10*} M \otimes \gg^{10} + T^{01*}  M \otimes \gg^{01}$, where, $ T^{10} M $,  $T^{01} M$ are  the holomorphic and anti-holomorphic distributions of $M$ and
$\gg^{10}$, $\gg^{01}$ are the subalgebras of $\gg^\bC$ (both isomorphic to $\gg$), which are  generated by the vectors of type $(1,0)$ and $(0,1)$.  
The $\gg$-valued $1$-forms $A_x$, $x \in \cU$,  being real,  the potential $A$ is uniquely determined by 
the  {\it $(1,0)$-potential}   $A^{10}\colon \cU \to T^{10*} M \otimes \gg^{10}$. \par
\subsection{Instantons on quaternionic K\"ahler  manifolds} \label{definitioninstantons}
Let $(M, g, Q)$ be a  qk manifold and consider  the isomorphisms  $T^\bC_x  M \simeq \sH_x \otimes \sE_x$    described  in \S \ref{notation}.  The  space of complex $2$-forms $\L^2 T_x^{*\bC} M$ splits  into three irreducible $\SL_2(\bC){\cdot} \Sp_n(\bC)$ moduli:
 \beq \label{decomp} \L^2 T_x^{*\bC} M  \simeq (\bC  \o_{\sH_x}) \otimes S^2 \sE^*_x + S^2 \sH^*_x \otimes (\bC  \o_{\sE_x})  + S^2 \sH^*_x \otimes \L^2_0 \sE^*_x\ .\eeq
Here $\o_{\sH_x}$ and $\o_{\sE_x}$ are the $\SL_2(\bC)$ and $\Sp_n(\bC)$  invariant symplectic forms of $\sH_x$ and $\sE_x$, respectively, and   $ \L^2_0 \sE_x$ is the irreducible $\Sp_n(\bC)$-submodule of $\L^2 \sE_x$  that is  complementary to $\bC \o_{\mathbf E_x}$.  The decomposition \eqref{decomp} is independent of the considered isomorphism $T^\bC_x  M \simeq \sH_x \otimes \sE_x$, since the latter
is unique up to the action of an  element in $\Sp_1{\cdot} \Sp_{p,q}$. 
\par
The ($\bC$-linear extension of the)  curvature tensor $F_x$   of  a gauge field  $(E, D)$  decomposes according to \eqref{decomp}. We write 
 \begin{align} \label{2.7}
 F_x = F^{(1)}_x + F^{(2)}_x\quad
  \text{with}\quad &F_x^{(1)} \in \o_{\sH_x}{\otimes} S^2 \sE^*_x \otimes \End(E_x)  
\nonumber  \\
 & F_x^{(2)} \in (S^2 \sH^*_x {\otimes}(\bC \o_{\sE_x}+ \L^2_0 \sE^*_x))\otimes \End(E_x)\, . 
  \end{align}
 A gauge field  $(E, D)$ on a qk manifold is called an {\it instanton} if the $F^{(2)}$ component  of its curvature tensor  vanishes identically.
Such instantons provide minima of the Yang-Mills functional $ \int_M |F|^2 \o_g$ and are thus Yang-Mills fields.  Instanton equations  on quaternionic K\"ahler manifolds have been studied by several authors, e.g.\  \cite{ward84,cgk,Ni,MS,Ti}, and are examples of a larger class of algebraic equations for the curvature tensor which imply the Yang-Mills equations \cite{cdfn, ACD, De}. \par
If an  instanton $(E, D)$ has  a compact (resp. complex) structure  group, we call it an {\it  instanton of compact} (resp. {\it complex}) {\it type}.
\par \medskip
\section{The Swann  bundle  of a quaternionic K\"ahler manifold} 
\label{swann}
\subsection{Conformally adapted frames}
Let $(M, g, Q)$  be a $4n$-dimensional qk manifold. Note that  each  homothetic Riemannian manifold 
$(M, r g)$  determined by a homothetic parameter  $r \in (0, + \infty)$  is  also a  qk manifold   
provided  it is associated  with the same bundle  of  quaternionic structures $Q \subset \End(TM)$ as 
$(M, g, Q)$.  We define the    bundle  of  {\it conformally adapted frames} of $(M, g)$, denoted 
$ \cC O_g(M, Q)$, as the union of all  adapted frame bundles  of  all its associated homothetic  manifolds,
$$p\colon \cC O_g(M, Q) \= \bigcup_{r \in (0, + \infty)}  O_{r g}(M, Q) \longrightarrow M\ . $$
This is a principal bundle with    structure group   
$$\bR_{> 0} \times \Sp_1 {\cdot} \Sp_{p,q} =  \bR^* \times_{\bZ_2} (\Sp_1 \times_{\bZ_2} \Sp_{p,q}) \ .$$
For each homothetic metric $r g$, the canonical embedding   of the bundle of adapted frames into the bundle of orthonormal frames $O_{r  g}(M)$,
 $\imath_r\colon O_{r g}(M, Q) \hookrightarrow O_{r  g}(M)$, 
 can be used to  pull-back the  tautological  $1$-form $\Q^{(r)}$ and the Levi-Civita connection  form  
 $\O^{(r)LC}$  of  $O_{rg}(M)$  onto $O_{r g}(M, Q)$. We thus have two  canonically defined $1$-forms,
$$\q^{(r)} \= \imath_r^* \Q^{(r)}\ \quad\text{and}\quad  \o^{(r)} \= \imath^*_r \O^{(r)LC}\ ,  $$
 the first being   $\bR^n$-valued and the second  $\sp_1 {+} \sp_n$-valued.
Although these  $1$-forms clearly depend on the homothetic parameter $r$, they combine  nonetheless to determine
two  smooth  $1$-forms $\q$  and $\o$ on  $\cC O_g(M, Q)$.   We call them the  {\it tautological $1$-form} and the  {\it Levi-Civita
connection $1$-form} of $\cC O_g(M, Q)$.  \goodbreak
 Note that:
 \begin{itemize}
 \item[--] $\q$ is the natural tautological $1$-form of $\cC O_g(M, Q)$  provided that this bundle is  considered to be a subbundle of the bundle of linear frames $L(M)$ of $M$; 
 \item[--] $\o$ is a connection $1$-form on $\cC O_g(M, Q)$ provided that  $\sp_1 + \sp_n$ is considered
 to be  a subalgebra of $(\sp_1 + \sp_n) + \bR$ (i.e. of   the Lie algebra of the structure group of the bundle of conformally adapted frames).
 \end{itemize}
 \par
\subsection{The Swann bundle and its canonical structure of a hyperk\"ahler cone}
\label{Swannsection}
  The {\it  Swann bundle} of a qk manifold  $(M, g, Q)$ is the principal  bundle 
 $$ \pi\colon \cS  \=  \cC O_g(M, Q)/\Sp_{p,q} \to M$$
  with  structure group  $ \bH^*\!{/}\bZ_2 \ \left(\ \simeq \bR^* \times_{\bZ_2} (\Sp_1 \times_{\bZ_2} \Sp_{p,q})/\Sp_{p,q} \right) $. \par
  \smallskip
As we will shortly see, this bundle admits a very special  hypercomplex structure $(J_1$, $J_2$, $J_3)$ and  a distinguished 
pseudo-Riemannian metric $h$ which make the tuple  $(\cS, h, J_1, J_2, J_3)$  an hk manifold canonically associated with $(M, g, Q)$.
This hk manifold  plays a crucial role in our discussion.  The explicit construction of the quadruple of tensor fields $(h,J_1,J_2,J_3)$ is as follows.
 \par
 \smallskip
 For  each point $ (x, [u]) \in \cS =  \cC O_g(M, Q)|_x/\Sp_{p,q}\,,\ x \in M$, consider the $\bH$-valued $1$-form  
 \beq \label{omegaS} \o^\cS_{(x,[u])} \= \left(\sigma^*(\o)_{(x,[u])}\right)^{\bH}\! :\ T_{(x,q)} \cS \longrightarrow  \bH\ , \quad \text{where }\eeq
 \begin{itemize}
\item[a)]  $\s\colon \cU\subset \cS  \to \cC O_g(M, Q)$ is a  local section  of $\cS$ defined   around $(x, [u])$, 
\item[b)]  $(\cdot)^{\bH}\colon \bR +  \sp_1 + \sp_{p,q} \to \bH$  denotes the  canonical   projection of $\bR +  \sp_1 + \sp_{p,q}$  onto $\bR + \sp_1 \simeq \bH$.
 \end{itemize} 
The mapping \eqref{omegaS} is independent of the choice of $\s$. Further, the  collection of $1$-forms 
$\o^\cS_{(x, [u])}$  determines   a smooth  $\bH$-valued $1$-form $\o^\cS$, which is in fact  a connection form on the $\bH^*\!{/}\bZ_2$-bundle $\cS$. \par
The $\bH^*\!{/}\bZ_2$-invariant horizontal distribution $\cD \subset T \cS$,  given by the kernel of the 
$1$-form $\o^\cS$, is naturally equipped with a quadruple of tensor fields  $g^\cD, J^\cD_\a$, $\a = 1, 2, 3$, 
defined pointwise as follows.
For each $(x, [u]) \in \cS|_x$, $x \in M$, we introduce: 
\begin{itemize}
\item[--] the tensor field in $\cD^*_{(x,[u])} \otimes \cD^*_{(x,[u])}$ given by  
\beq\label{g} g^\cD_{(x,[u])} (v,w)  \= r g_x(\pi_*(v), \pi_*(w))  \eeq
\item[--]\label{J}  the tensor field in $\cD_{(x,[u])} \otimes \cD^*_{(x,[u])}$ given by 
\beq
J^\cD_\a|_{(x,[u])}(v)  \= (\pi_*|_{\cD})^{-1}(J^{(u)}_\a)\ , \label{theJ} \eeq
\end{itemize}
where  $(J^{(u)}_\a)$  is the unique triple of complex structures  of $T_xM$ satisfying \eqref{adapted} 
for a frame $u = (e_i)$  in equivalence class $[u]$. It is possible to check that such a  triple   does not depend on the 
particular adapted frame  $u = (e_i)$ chosen out of  the equivalence class $[u]$,  meaning that  \eqref{theJ} is well defined.
Using the above construction we may verify the following (see  e.g. \cite{Ha, PPS} for details): 
\begin{itemize}
\item[a)]   the   action of the subgroup  $\{e^t, t \in \bR\}$ of  $ \bH^*\!{/}\bZ$ leaves invariant each of 
the tensor fields $J_\a$, but not the tensor field $g^\cD$, which transforms as 
\beq R_{e^t}^* (g^\cD)= e^t\, g^\cD\ .\eeq
\item[b)]  the right action $R_{[q]}\colon \cS \to \cS$ of each element  of the structure group $\bH^*\!{/}\bZ_2$ of the form
$$[q] = q\ \text{mod} \ \bZ_2  \in \Sp_1{/}\bZ_2 \subset \bH^*\!{/}\bZ_2$$  
leaves invariant $g^\cD$, but not the triple  $(J_1^\cD, J_2^\cD, J_3^\cD)$, which get rotated  in the following sense:  at each point 
$(x,[u]) \in \cS|_x$, the triple $\left(R_{[q]*} (J^\cD_\a|_{(x, [u])})\right)_{\a = 1, 2, 3}$ is the  triple of complex structures on $\cD_{(x, [u]{\cdot} [q])}$ projecting onto the triple of complex structures of $T_x M$ 
   \beq \label{permuting}J^{(q)}_\a|_x \= R^{(q^{-1})}\circ J_\a|_x \circ R^{(q)} \ , \eeq
  where   $R^{(q)}$ is the  right $\Sp_1$-action
  $$R^{(q)}\colon  T_x M   \longrightarrow T_x M\ ,\qquad R^{(q)}(v) \=v{\cdot} q\ ,$$
determined by one of the  isomorphisms   $T_x M \simeq \bH^n$  established by an adapted frame 
$u = (e_i)$ in the equivalence class   $[u]$. 
It can be directly checked that the complex structures \eqref{permuting} are independent of the choices of $q \in [q]$ and  $u \in [u]$.
\end{itemize}
We now   canonically  ``extend'' the tensor fields $g^\cD$ and $J_\a^\cD$, which we have so far defined only 
on vectors of  the distribution $\cD \subset T\cS$, to all of  $T \cS$.  
First recall  that each $1$-form $\o^\cS|_{(x, [u])}$   gives an isomorphism between  the vertical  tangent space 
$\cV_{(x,[u])} = T^V_{(x,[u])} \cS $ and    $\bH = \Lie(\bH^*\!{/}\bZ_2)$. It therefore  induces  a natural   hypercomplex
structure $(J_\a^\cV|_{(x,[u])})$ on  $\cV_{(x,[u])}$  corresponding to the standard hypercomplex  structure of 
$\bH$, i.e. $(J_1^\cV|_{(x,[u])} \= {\bf i}, J^\cV_2|_{(x,[u])}\=  {\bf j}, J^\cV_3|_{(x,[u])} \= {\bf k})$. 
We  thus have a triple  $(\bJ_\a)_{\a = 1,2,3}$ of almost complex structures  on $\cS$ defined at 
each point $(x, [u])$  by 
\beq \label{triple1} \bJ_\a|_{(x,[u])}= J^\cD_\a|_{(x,[u])} + J^\cV_\a|_{(x,[u])}\ .\eeq
The bundle $\cU_{n/(n+1)}(M) \to M$,    considered by Pedersen, Poon and Swann  in \cite{PPS}  (or,  equivalently,    the bundle ${\mathbf H}' \to M$  considered in  \cite{Sa1})  is identifiable with the 
Swann bundle $\cS \to M$ defined  above and  the triple  \eqref{triple1} coincides  with the triple of  almost complex structures  defined there.  Thus,  by  Prop. 3.4 of \cite{PPS}, we immediately have  the following 
 \begin{prop}[\cite{PPS}] \label{prop13} The   almost complex structures  \eqref{triple1} are   integrable and  
 form  a hypercomplex structure   on $\cS$  satisfying the following  conditions: 
 \begin{itemize}
  \item[1)] The right action of the subgroup $\{e^t, t \in \bR\}\subset \bH^*\!{/}\bZ_2$  leaves each $\bJ_\a$
  invariant, whereas the right action of $\Sp_1{/}\bZ_2 \subset \bH^*\!{/}\bZ_2$ {\rm rotates} them according to 
  \eqref{permuting}.
 \item[2)]  The distributions $\cV$,  $\cD$  are   invariant under each complex structure  $\bJ_\a$.
 \item[3)] Along every  fibre $\cS|_x  \simeq \bH^*\!{/}\bZ_2$, the triple $(\bJ_1, \bJ_2, \bJ_3)|_{\cS_x}$  
 is  equal to  $({\bf i}, {\bf j}, {\bf k})$. 
 \item[4)] For each $(x, [u]) \in \cS|_x$, the projection $\pi_*$ maps   the  triple $(\bJ_\a|_{(x,[u])})$  onto 
 a triple $(J_\a)$ of  complex structures of $T_x M$ belonging to the twistor bundle  $\ZZM$ of $(M, g, Q)$,
 thus establishing  a   $\bC^*$-bundle  $p\colon \cS \to \ZZM$  over $\ZZM$.
   \end{itemize}
 Up to a bundle automorphism  preserving  $\o^\cS$, the triple $(\bJ_1, \bJ_2, \bJ_3)$ is 
 the unique hypercomplex structure on $\cS$ satisfying 1) -- 4). 
 \end{prop}
 \par
 \medskip
Having  ``extended"  the triple of tensor fields $J^\cD_\a$ in a canonical way, we  now similarly  ``extend''  $g^\cD$.
 Let $\langle \cdot, \cdot\rangle$ be the standard $\Sp_1$-invariant Euclidean  product on $\bR^4 \simeq \bH$ and 
on each vertical subspace $\cV_{(x,[u])} \subset T_{(x,[u])} \cS$, define the inner product  
$$g^\cV(v, w)_{(x,[u])} \= \langle \o^\cS(v), \o^\cS(w) \rangle\ ,\qquad v, w \in T^V_{[u]} \cS\ .$$
If  $(M, g, Q)$ is a strict qk manifold, i.e. with scalar curvature $\L \neq 0$, the {\it canonical Swann metric}  
on  $\cS$ is the  pseudo-Riemannian metric  
defined by 
\beq\label{swannmetric}
\begin{split}
h(X, Y) =
 \left\{ 
 \begin{array}{ll} \frac{1}{\Lambda}g^\cV(X, Y) & \text{if}\ X, Y \in\cV\ ,\\[6pt]
0 & \text{if}\ X \in\cV,\ Y \in \cD\ ,\\[6pt]
g^\cD(X, Y) & \text{if}\ X, Y \in \cD\ .
\end{array}
\right. 
\end{split}
\eeq
If $(M, g, Q)$ is non-strict ($\L = 0$), the {\it canonical Swann metric}   on $\cS$ is taken to be any  
metric  having the form   \eqref{swannmetric}  with  $\frac{1}{\L}$ replaced by  some  positive constant 
$c > 0$.  For a  canonical Swann metric $h$ the following hold:
\begin{itemize}[itemsep=5pt plus 4pt minus 2pt]
\item[i)] {\it $h$  is Hermitian with  respect to each of the three complex structures $\bJ_\a$}. 
\item[ii)] {\it If $(M,g)$ has  signature $(4p, 4q)$, the signature of $h$ is $(4p+4, 4q)$ if $\L \geq 0$ 
and $(4p, 4q + 4)$ if $\L < 0$}.
\item[iii)] {\it $h$ is invariant under the right action of every  $q \in \Sp_1{/}\bZ_2\subset \bH^*\!{/}\bZ_2$}. 
\item[iv)] {\it  For each  $e^t \in \bH^*\!{/}\bZ_2$, we have that  
$R_{e^t}^* h = e^t\, h$} (i.e.  $h$ behaves  as a  conical metric with respect to the dilatations $e^t$).
\item[v)] {\it The quadruple $(h, \bJ_1, \bJ_2, \bJ_3)$ is an hk structure on $\cS$.} 
\end{itemize}
Note that if   $\L = 0$,  property (v)  is essentially a trivial consequence of the fact that  $\cS$ is a locally trivial bundle over $M$. If  $\L \neq 0$, (v) is a non trivial property, following from results  in  \cite[\S 3]{Sw}.\par
\smallskip
The  hk manifold $(\cS, h, \bJ_1, \bJ_2, \bJ_3)$  is  called the {\it hyperk\"ahler (hk)  cone} of the qk manifold  $(M, g, Q)$.\par
 \subsection{$\bH^*$-conical manifolds and their relation to qk manifolds}
 \label{hkcone}
We have seen in the previous section that any (strict) qk manifold has a canonically associated hk manifold, namely its hk cone $(\cS, h,  \bJ_\a)$. 
This correspondence has  an  inverse. Let us first consider the following:
\begin{definition} \label{conicalH}  Let   $(N, h,  \bJ_\a)$ be  an hk manifold. We say that it is     
{\it $\bH^*$-conical}  if  it is equipped with an almost effective right  action 
$R\colon \bH^* \to \operatorname{Diff}(N)$,  with $\ker R = \bZ_2$,  and displays  the following  properties.
\begin{itemize} [itemsep=8pt plus 5pt minus 2pt]
\item[a)] The effective  right  $ \bH^*\!{/}\bZ_2$-action makes $N$ an $ \bH^*\!{/}\bZ_2$-bundle over \linebreak $M\, {=}\, N{/}\bH^*$
such that the  restrictions of $h$  and  $\bJ_\a$ to the vertical tangent spaces render  each fibre  $N_x$,
being $3$-holomor\-phi\-cally isometric to $(\bH^*\!{/}\bZ_2\,\,,  c \langle\cdot\,, \cdot \rangle)$ for some constant $c$,
independent of $x$.  
\item[b)]  The  distribution $\cD$, given by the spaces that are $h$-orthogonal to  the $\bH^*$-orbits in $N$,
is the horizontal distribution of a connection form on the bundle $N \to M {=} N{/}\bH^*$.
\item[c)]  For   each   $\{e^t,\ t \in \bR\} \subset \bH^*$, the map $R_{e^t} \colon N \to N$  is $\bJ_\a$-holomorphic
for each $\bJ_\a$,  but it acts on  the metric $h$ as a homothety, $R_{e^t*} h = e^t h$.
\item[d)] For each $[q]  \in  \bH^*\!{/}\bZ_2$ of the form $[q] = q\ \text{mod}\ \bZ_2$ for  $q \in \Sp_1{/}\bZ_2$, 
the map $R_{q} \colon N \to N$ acts on $h$ and  $\bJ_\a$ as follows: 
\begin{itemize} 
\item[--] it rotates the hypercomplex structure $(\bJ_1, \bJ_2, \bJ_3)$  as in  \eqref{permuting}
\item[--] it leaves $h|_{\cD \times \cD}$ invariant.
\end{itemize} 
\item[e)] If  $H_D$ is the vector field on $N$ generating the $1$-parameter family of diffeomorphisms 
$\Phi^{H_D}_t \= R_{e^t}$, the tensor field $\gc \in \L^2 \cD^* $ defined by 
\beq\label{tensorc}\ \hskip 1.2cm \gc_y(v, w) \=  h_y(\n_v H_D, w) -  h_y(\n_w H_D, v)\ ,\quad y \in N,\ \ v,w \in \cD_y,\eeq
vanishes identically.  
\end{itemize}
\end{definition}
Note   that   c) implies  the invariance of the Levi-Civita connection $\n$ under  diffeomorphisms $R_{e^t}$.  
Thus  {\it condition  e)  is equivalent to requiring that $\gc$  vanishes  identically on  at least one  
hypersurface  transversal to the orbits of  $\{e^{it}, t \in \bR\}$}. \par
\medskip
The following theorem gives a characterisation of hk cones. It contains  a  reformulation of  certain results in  \cite{PPS, Sw, BCDGVV}. \par
\begin{theo} \label{theorem25} Let $(N, h,  \bJ_\a)$ be a $4(n+1)$-dimensional $\bH^*$-conical 
manifold and  $\cV^{\Sp_1}$ and $\cD $  the distributions in $TN$ consisting respectively of vectors tangent to
the  $\Sp_1$-orbits and those  $h$-orthogonal to the $\bH^*$-orbits.  Then:
\begin{itemize}
\item[i)] The  distribution  $\cD + \cV^{\Sp_1}$ is involutive.
\item[ii)]  Given a maximal integral leaf $S \subset N$ of $\cV^{\Sp_1} + \cD$ and denoting 
  $M \= S/\Sp_1$, define  the pseudo-Riemannian metric $g$ and the bundle $Q \ \subset \End (TM)$  
  over $M$    by  
\beq \label{metric} \begin{split} &g_x(v, w) \=h_{y}\left(v,  w\right)\ ,\\
 &Q_x \= \operatorname{span}_\bR \{ J = \bJ_\a|_{\cD_{y}}\ ,\ y \in x =  [y]\ \} \ ,
\end{split}\eeq
for each $x = [y] \in M = S/\Sp_1$ (using of course the isomorphism $\pi_*\colon TS \to T M$ to identify each $\cD_y \subset T_y S$ with $T_{x=[y]}M$). 
Then, $(M, g, Q)$ is a qk manifold, whose  hk cone is  3-holomorphically isometric to  $(N, h, \bJ_\a)$. 
\end{itemize}
\end{theo}
\par 
 \subsection{A geometrical interpretation of the Swann bundle}
\label{geometricinter}
 As discussed in \S \ref{notation}, the adapted frames of a qk manifold $(M, g, Q)$ establish  $\bC$-linear  isomorphisms between the complexified tangent spaces  $T^\bC_x M$ and $\bC^2 \otimes \bC^{2n}$  and make  $T^\bC_x M $  an $ (\SL_2(\bC) \times_{\bZ_2} \Sp_n(\bC))$-module  of the form $T^\bC_x M \simeq \sH_x \otimes \sE_x$, with    
 $$\sH_x = T^\bC_x M/\SL_{2n}(\bC) \simeq \bC^2\ ,\qquad \sE_x = T^\bC_x M/\SL_{2}(\bC) \simeq   \bC^{2n} .
 $$ 
We recall that  the  isomorphism  $T^\bC_x M \simeq \sH_x \otimes \sE_x$ is independent  of the 
adapted frame $u = (e_i)$  used in its construction.  Thus, the family  of  vector spaces 
$\sH_x \simeq \bC^2$, $x \in M$,   constitutes a rank two holomorphic vector bundle $\pi\colon \sH \to M$, intrinsically determined by  the quaternionic structure of $(M, g, Q)$ and locally  
equipped with a field $\o_{\sH}$ of  $\Sp_1$-invariant  symplectic $2$-forms along the fibres (see \cite{Sa}). The field $\o_{\sH}$ is
 uniquely determined up to multiplication by a nowhere vanishing function.\par

Now,  let us denote by $(h^o_1, h^o_2)$, $(e^o_a)$ and $(e_{ia}^o)$ the standard bases for 
$\bC^2$,  $\bC^{2n}$ and $ \bC^2 \otimes \bC^{2n}$, respectively, i.e.  let  
$$\big(h^o_1 \= (1,0), h^o_2 \= (0,1)\big)\ ,\ \ 
\big(e^o_a \= (0, \ldots, \underset{\text{$a$-th entry}}1, \ldots, 0)\big)\ ,\ \  
\big(e_{ia}^o \= h^o_i \otimes e^o_a\big).
$$ 
Further,  for any ($\bC$-linearly extended) adapted frame $u\colon  (\bH^n)^\bC = \bC^2 \otimes \bC^{2n}  \to T^\bC_x M$,  denote  the corresponding bases for   $T^\bC_x M$, $\sH_x$ and $\sE_x$ by 
$$( e_{ia}  \= u(e^o_{ia})) \ ,\quad (e_i \= u(e^o_{ia}) \!\!\mod \SL_{2n}(\bC)) \ ,\quad (e_a \= u(e^o_{ia})\!\!\mod  \SL_{2}(\bC)). 
$$

Note that when an adapted  frame  $u\colon \bH^n \to T_x M$ is changed into   another adapted   frame, 
 $$u' = u \circ U\qquad \text{with}\quad U = \left(\begin{array}{cc} u^1_+ & u^1_-\\ u^2_+ & u^2_-\end{array} \right)  \in \Sp_1\subset  \SL_2(\bC) = \Sp_1(\bC)\ ,$$
 the  corresponding complex basis $(e_{ja})$ of $T^\bC_x M$ changes as
 \beq \label{25} e_{+a} = u^j_+ e_{ja}\ ,\qquad e_{-a} = u^j_- e_{ja}\ , \eeq
and  the corresponding basis $(e_i)$ for the $2$-dimensional space $\sH_x  \simeq \bC^2$  transforms into
 \beq \label{25bis} e_{+} \= u^i_+ e_{i}\ ,\qquad e_{-} \= u^i_- e_{i}\ .\eeq
 We shall frequently use the complex frames $(e_{\pm a})$ and $(e_\pm)$ and for brevity,  we 
 call them the  {\it  adapted complex frames} for  $T^\bC_x M$ and $\sH_x$, respectively. \par
 \smallskip
 Note that each adapted complex frame $(e_{+}, e_{-})$ for $\sH_x  \simeq \bC^2$ is orthonormal and 
 symplectic and is therefore uniquely determined, up to the  sign of the second vector $e_{-}$, by just 
 its first vector $e_+$. This means that the projection 
 \beq \label{projectj} j\colon \cC O_{g}(M, Q) \to \sH^o \= \sH {\setminus} \{\text{zero section}\}\ ,\ \  
 j(t e_{+a}) \=  t e_+\ ,\quad t \in (0, + \infty)\ ,
 \eeq
 determines a fibre preserving diffeomorphism between the Swann bundle $\pi\colon \cS \to M$ and the  holomorphic bundle $p\colon \sH^o \to M$,  
 \beq \label{tildej} \wt j\colon  \cS = \cC O_{g}(M, Q)/\Sp_{p,q} \to \sH^o\ .\eeq
 We may therefore  geometrically identify  $\cS$ with  the  bundle $\sH^o \to M$,  
 {\it provided that $\sH^o$ is equipped with the  structure of a principal $\bH^*\!{/}\bZ_2$-bundle 
 determined by  the  right action of  
 $ \bR^* \times_{\bZ_2} (\Sp_1 \times_{\bZ_2} \Sp_{p,q})$ on $\cC O_g(M, Q)$ through the projection  \eqref{projectj} of $\cC O_g(M, Q)$  onto $\sH^o$. }
 \par
 \subsection{$\bH P^n$ and its hyperk\"ahler cone}
 \label{Hpn}
 \label{section35}  Consider now the quaternionic projective space  $\bH P^n =  (\bH^{n+1} {\setminus} \{0\}){/}\bH^*$,   with $\bH^*$  acting on the right, and  the following  related notational data. 
 \begin{itemize} 
 \item[1)]  We denote by $\langle q, q'\rangle = \sum_{i =0}^n \overline q^i q'{}^i$   the standard Hermitian product of
   elements  $q, q' \in \bH^{n+1}$ and  by  $(e^o_j)$  the standard orthonormal  basis of $\bH^{n+1} \simeq \bR^{4(n +1)}$ with   
  \begin{align*}
 e^o_{4 k +1} &\= (0, \ldots, 0,  \underset{(4 k +1)\text{-th place}} 1, 0 \ldots , 0)\ ,\\
 e^o_{4 k +2} & \= (0, \ldots, 0,  \underset{(4 k +1)\text{-th place}} {\bf i}, 0 \ldots , 0)\ ,\\
 e^o_{4 k +3} &\= (0, \ldots, 0,  \underset{(4 k +1)\text{-th place}} {\bf j}, 0 \ldots , 0)\ ,\\
 e^o_{4 k +4} & \= (0, \ldots, 0,  \underset{(4 k +1)\text{-th place}} {\bf k}, 0 \ldots , 0)\ ,\qquad  \text{with }\ \ 0 \leq k  \leq n\ .
\end{align*}
 \item[2)] $S(\bH^{n+1})$ is the unit sphere of $\bH^{n+1}$.  The  natural projection   $\pi\colon S(\bH^{n+1}) \to \bH P^n \simeq   S(\bH^{n+1})/S(\bH) $ is the    well-known  {\it Hopf fibration}.
 \item[3)] $g_o$ is   the    Riemannian metric of $\bH P^n$   defined   by   projecting   each  inner product   
 $\langle \cdot, \cdot \rangle|_{\cD_q \times \cD_q}$  of  the spaces  of the distribution 
 $\cD  \subset TS(\bH^{n+1})$ orthogonal  to the  fibres  of the Hopf fibration  onto the corresponding  tangent space $T_{[q]} \bH P^n$,  $[q] \in \bH P^n$. 
 \item[4)]  $Q \subset \End(T \bH P^n)$ is the bundle of  endomorphisms of the tangent spaces determined as follows. For each $[q] \in \bH P^n$, consider the points $q \in S(\bH^{n+1}) \subset  \bH^{n+1} {\setminus} \{0\}$  projecting onto $[q]$ 
 and the triples of complex structures on the horizontal spaces $\cD_q\subset T_qS(\bH^{n+1})$
 determined by the right multiplications  by  ${\bf i} $, ${\bf j}$ and ${\bf k}$.  Each  triple projects onto a different hypercomplex structure, but each such projected hypercomplex structure spans a $3$-dimensional vector space $Q_{[q]} \subset \End( T_{[q]} \bH P^n)$,  which is  independent  of  $q \in [q]$.
 The bundle $Q$ is defined as    $Q \= \bigcup_{[q] \in \bH P^n} Q_{[q]}$. 
\end{itemize}
 The triple $(\bH P^n, g_o, Q)$ is a typical   example of  a homogeneous qk manifold.  We remark  that the standard isometric left action of $\Sp_{n+1}$ on $\bH^{n+1}$ acts transitively and effectively on each sphere $t S(\bH^{n+1})$ of radius $t > 0$  and with isotropy $\Sp_{n}$.
This transitive left action  commutes with the right action of $\Sp_1$ and gives rise to the
$\Sp_{n+1}$-equivariant map 
 $$j^{(t)}\colon t S(\bH^{n+1})  \simeq \Sp_{n+1}/\Sp_n \longrightarrow \bH P^n \simeq \Sp_{n+1}/\Sp_n\times_{\bZ_2} \Sp_1$$
 for each $t \in (0, + \infty)$.  
  The left  $\Sp_{n+1}$-action on $\bH P^n$ is isometric  for each homothetic metric $t g_o$ and  it  naturally lifts to the level of linear frames, mapping the bundle of  adapted frames   into itself. All this allows the identifications  
  $$O_{t g_o}(\bH^n,  Q) \simeq  \Sp_{n+1}{/}\bZ_2\ ,\qquad \cC O_{g_o}(\bH^n,  Q) \simeq   \bR_{> 0} \times \Sp_{n+1}{/}\bZ_2 $$
 and shows the existence of an  $\Sp_{n+1}$-equivariant diffeomorphism
\beq
j\colon \bH^{n+1} {{\setminus}}\! \{0\} {/}\bZ_2 = \!\! \bigcup_{t \in (0,+\infty)} {\hskip -4pt}  t S(\bH^{n+1}){/}\bZ_2 \longrightarrow  \cS 
= \cC O_{g_o}(\bH^n,  Q)/\Sp_n\, .
\eeq
This implies  that  the Swann bundle $\pi\colon \cS \to \bH P^n$  is equivalent to $$p\colon (\bH^{n+1}{\setminus}\{0\}) {/ }\bZ_2 \to \bH^n\,,$$
the $\bZ_2$-quotient of the tautological bundle of $\bH P^n$.
Note that both the flat  metric  $h_o = \langle \cdot , \cdot \rangle$
and the standard hypercomplex structure $( {\bf i},  {\bf j},  {\bf k})$ of $\bH^{n+1} {\setminus} \{0\}$, 
given by the right multiplications by the three complex structures of $\bH$, are invariant under the left 
$\Sp_{n+1}$-action. Using this, it can be  checked that the {\it right action} of $\bH^* = \bR_{> 0} \times \Sp_1$ on $(\bH^{n+1}{\setminus} \{0\}){/}\bZ_2$ transforms  $h_o$ and    $( {\bf i},  {\bf j},  {\bf k})$  as  prescribed  for  
$\bH^*$-conical manifolds, so that  $((\bH^{n+1} {\setminus} \{0\}){/}\bZ_2, h_o,    {\bf i},  {\bf j},  {\bf k})$   is identifiable  with   the   hk cone of $(\bH P^n, g_o, Q)$.

\par\medskip
\section{Harmonic spaces}
\setcounter{equation}{0}
\subsection{Harmonic space of a hyperk\"ahler manifold}
\label{sect41}
Let   $(N, g, J_1, J_2, J_3)$  be an hk manifold, i.e. a qk manifold  admitting  a $\n$-parallel hypercomplex 
structure  $(J_1, J_2, J_3)$  that  generates  the fibre $Q_x \subset Q$ at each $x \in N$.  Then,  for   any point  
$z = (a, b, c) \in \bR^3$ which is in  the unit sphere $S^2 \simeq \bC P^1$, the linear combination 
$ I^{(z)} \= a J_1 + b J_2 + c J_3$  is also a $\n$-parallel   integrable complex structure.  
The {\it twistor bundle}  $\ZZ$ of $(N, g, J_1, J_2, J_3)$ is    the trivial $\bC P^1$-bundle  over $N$ defined by 
$$\pi^\cZ\colon  \ZZ \=  N \times  \{I^{(z)}, z\in S^2\}  \simeq N \times \bC P^1 \to  N\ .$$
It is naturally equipped with an integrable complex structure  $\wh I$, coinciding with  the complex 
structures $I^{(z)}$ along  the horizontal leaves $N \times \{z\}$, as well as with the standard $\SL_2(\bC)$-invariant
complex structure of $\bC P^1$ along the vertical spaces  $\{x\} \times \bC P^1 $ \cite{Sa, Hi}.
\par
The {\it harmonic space} $\HN$ of $(N, g, J_1, J_2, J_3)$
is an analogous trivial bundle, 
albeit with  larger fibres.  It  is defined as  the trivial $\SL_2(\bC)$-bundle   
\beq\label{H_hk}
\pi^\cH\colon \HN \= N\times \SL_2(\bC) \to N\ ,
\eeq
equipped with the  following  integrable complex structure. 
For each  $(x, U) \in \HN$, consider  the  direct sum decomposition
$T_{(x, U)} \HN =  T_x N + T_U \SL_2(\bC) \simeq  T_x N +  \sl_2(\bC)$,  where  $T_U \SL_2(\bC)$
is identified with  $\sl_2(\bC)$ by using  {\it right invariant}  vector fields (\footnote{More precisely,  
we assume that any  $v \in T_U \SL_2(\bC)$ is identified with the unique element 
$E^{(v)} \in \sl_2(\bC)$, whose associated  right invariant vector field $\wh{E^{(v)}}$ is such that  
$\wh{E^{(v)}}|_U = v$.}). 
Now,  denote by $\bI_{(x,U)}$   the unique linear map on  $T_{(x, U)} \HN$ satisfying
\beq
\bI_{(x,U)}|_{T_x N} \= I^{(z)}|_x,\  \ \ z = U{\cdot}[0:1] \in \bC P^1\ ,  \qquad \bI_{(x, U)}|_{\sl_2(\bC)} = \Jst\ ,
\label{sl2cxstr}\eeq
where  $\Jst$ is  the  complex structure of $\sl_2(\bC)$ given  by the  multiplication by 
$\left(\smallmatrix i & 0\\ 0 & i\endsmallmatrix	\right)$. The collection of linear maps  $\bI_{(x,U)}$ defines an 
almost complex structure $\bI$ on $\HN$, which can be checked to be integrable. 
The complex manifold   $(\HN , \bI)$   naturally projects  to  the twistor space $(\ZZ, \wh I)$,   the  projection 
being holomorphic with fibres given by the  orbits in $\HN$ of the subgroup  $B \subset \SL_2(\bC)$  of 
upper triangular matrices (see  e.g. \cite[\S 3.2]{DPS}).\par

\subsection{Harmonic space of a quaternionic K\"ahler manifold}

We now proceed to the harmonic space of a strict qk manifold.
One might expect this to be defined analogously to the hyperk\"ahler case: 
The harmonic space $\mathcal H (N)$ of an hk manifold $N$ is a topologically trivial bundle 
equipped  with a non-product complex structure which makes it  a holomorphic bundle over its 
twistor space $Z(N)$. Now, the twistor construction for any qk manifold was introduced by 
Salamon \cite{Sa} regardless of the sign of its scalar curvature as a $\bC P^1$-bundle in 
analogy with the Penrose construction in real dimension four. As a matter of fact it applies 
equally well  in the more general context of quaternionic manifolds, see \cite{AMP} for example.
However, this twistor space does not serve as a suitable base for a prospective 
harmonic space for a strict qk manifold $M$, having the drawback that the lifted gauge field  is not 
partially flat in the sense required for the application of the  harmonic space method 
(see e.g. \cite{gio}). This difficulty can be circumvented by considering instead of the qk manifold 
its hyperk\"ahler cone.
We thus have:
\par
\begin{definition} 
The  {\it harmonic space  of a (strict) qk manifold $(M, g, Q)$} is the  bundle over $M$  given by the pair  $( \HMqk,  \pi^{\Hqk})$,  where 
 \begin{itemize}
 \item[(i)] the total space  $\HMqk \= \cH(\cS) $  is the harmonic space of  the  hk cone 
 $(\cS, h, \bJ_\a)$  of $M$ (cf. \eqref{H_hk}) and 
 \item[(ii)]  the projection $ \pi^{\operatorname{qk}}$ is  the map 
 $$ \pi^{\operatorname{qk}} = \pi \circ \pi^{\cH|\cS}:  \HMqk = \cS \times \SL_2(\bC) \to M \,,  $$ 
 determined by  composing  the  projection  $\pi^{\cH|\cS}: \HMqk = \cH(\cS) \to \cS$ of the harmonic space $\cH(\cS)$ onto $\cS$ with the projection  $\pi: \cS \to M$ of the 
 Swann bundle over $M$.
 \end{itemize}
 In other words, $\HMqk$ is precisely the harmonic space  $\HMqk = \cH(\cS)$ of $\cS$,  
 considered as a bundle over  $M$ rather than one over $\cS$. 
 \end{definition} 

 The following lemma gives a few basic properties of $\HM^{\operatorname{qk}}$.
\begin{lem} The manifold $\HMqk$  has  a natural   integrable complex structure
and is acted  on holomorphically  and freely by the group 
\beq \label{definitionL} L \=  (\bH^*/\bZ_2) \ltimes \SL_2(\bC)\ . \eeq
The action is transitive on  the fibres of the projection   
$\pi^{\operatorname{qk}}\colon \HMqk \to M$ and gives $ \HMqk$
the structure of a  principal  $L$-bundle over $M$.
 \end{lem}
 \begin{pf}
Since $\HMqk$ is  the harmonic space of an hk manifold,  the discussion  in  \cite[\S 3.2]{DPS}  implies  that it  is  naturally  equipped with an integrable complex structure with respect to which  the orbits of the right action of $\SL_2(\bC)$ are complex. Further, $\HMqk$
can be (locally) identified  with the bundle of adapted complex  frames $(e_+, e_-)$ for the rank $2$ holomorphic 
vector bundle $\sH \to \cS$ defined  in \S \ref{geometricinter}.  
Now,  for  any element $a \in \bH^*/\bZ_2$, let us  consider the corresponding diffeomorphism  $f_a\colon \cS \to \cS$
determined by the right action of $\bH^*/\bZ_2$ on $\cS$  described in detail in \S \ref{Swannsection}. From this
description it follows immediately  that the  push-forward  $f_a{}_*\colon T \cS \to T\cS$
transforms any complex adapted frame  of $(\cS, h, J_\a)$ into another complex adapted frame. 
This means  that the push-forward map $f_a{}_*$  determines  an automorphism
of the bundle of the vertical frames $(e_+, e_-)$ of $\sH$ which preserves its complex structure.
This is tantamount to saying  that each element  $a \in \bH^*/\bZ_2$ is uniquely associated with a  fibre preserving  biholomorphism of
$\pi^{\cH|\cS}\colon \HMqk \to \cS$. 
These biholomorphisms  together with  those of the  $\SL_2(\bC)$-action determine an effective right action of
$L = \SL_2(\bC) {\ltimes} (\bH^*/\bZ_2)$ on  $\HMqk$  acting transitively on the fibres of $\HMqk$ over $M$.
\end{pf}
\begin{rem}
We recall that  $\HMqk$, although equipped with a non-product complex structure, is topologically a cartesian
product of $\SL_2(\bC)$ and $\cS$. Now,  $\cS$ is an $\bH^*\!/\bZ_2$-bundle 
over $M$. So the  product group   $ \SL_2(\bC) {\times} (\bH^*/\bZ_2)$  acts naturally and 
transitively on the fibres of $\pi^{\text{qk}}\colon \HMqk \to M$. 
This natural action  restricts to an action of the  $\bH^*/\bZ_2$ subgroup on $ \HMqk$, which leaves 
each horizontal leaf $\cS {\times} \{U\} \subset  \cS {\times} \SL_2(\bC) = \HMqk$ invariant. 
However, this action of  $\bH^*/\bZ_2$ is not  holomorphic. In contrast,  the action of the semi-direct product 
$L = (\bH^*/\bZ_2) \ltimes \SL_2(\bC)$ described in the  above proof  is holomorphic and hence restricts to 
a new  action of its subgroup $\bH^*/\bZ_2$ which is also holomorphic. 
The price we pay is that  this action   does not  leave  the  horizontal leaves $\cS {\times} \{U\} $ invariant.  
In fact each element $a \in  \bH^*/\bZ_2$ biholomorphically maps each leaf $\cS {\times} \{U\} $  to  another leaf 
$\cS {\times} \{U'\} $,   with $\,U' = U{\cdot} B_a\,$ for some element $B_a \in \SU_2$. 
 \end{rem} 
 
From the discussion in  \S \ref{section35}, the  harmonic space $\HMqk$ of $M = \bH P^n$ is identifiable with the  bundle    
$$ \pi^{\operatorname{qk}} \colon (\bH^{n+1} {\setminus} \{0\}/\bZ_2) \times \SL_2(\bC) \to \bH P^n = (\bH^{n+1} {\setminus} \{0\})/\bH^* .$$
Note that  this $\HMqk$  can be seen as a principal bundle with respect to both the  action of the direct product 
$(\bH^{n+1} {\setminus} \{0\}/\bZ_2) \times \SL_2(\bC))$, as well as the action of the semi-direct product 
$ L = (\bH^{n+1} {\setminus} \{0\}/\bZ_2) \ltimes \SL_2(\bC))$, depending on  whether $\HMqk$ is considered as a cartesian product 
or as a quotient of  an appropriate  bundle of linear frames.   
\par
\subsection{Complexified harmonic space} \label{cxH}
Consider an  $n$-dimensional complex manifold   $(N, J)$. 
The {\it complexification} of $(N, J)$  is the pair $(N^{\bC}, \imath)$ given by:
 \begin{itemize}
\item[a)]  the complex manifold  $N^{\bC} \= N \times N$ having  complex  structure $\wt J$  defined at each point 
$(x, y) \in N{\times} N$ by $\wt J_{(x,y)} (v, w) {\=} J_x(v) {-} J_y(w)$
\item[b)] the standard  diagonal embedding 
$\imath\colon N \to N^{\bC}$,  $\imath(x) = (x, x)$.
\end{itemize}
Note that   $\wt J$ is so defined that the complete atlas of  holomorphic coordinates of $(N^\bC, \wt J)$ is
precisely the one generated by the collection of complex coordinates  of  the form  
$\wt \xi = (z^i, \overline{z'{}^j})\colon \cU \times \cV \to \bC^{2n}$  determined by pairs of holomorphic coordinates   
$(z^i)$,  $(z'{}^j)$  of $(N, J)$.
 \par
\smallskip
This notion  of   ``complexification''  of  a complex manifold  yields   
the following notion of ``complexified harmonic space''.
Consider the harmonic space $\HMqk = \cS \times \SL_2(\bC)$   of a (strict) qk manifold $M$ and denote by  $\HCMqk$
 the cartesian product  $\HCM^{\operatorname{qk}} = \cS \times \cS \times \SL_2(\bC)$  equipped with the unique (integrable) complex structure 
 $\bI^{(\bC)}$  which coincides with the right invariant  complex structure  along the  (vertical) leaves
$\{x\} \times\{y\} \times \SL_2(\bC) \simeq \SL_2(\bC)$  (see  \eqref{sl2cxstr}) and with the complex structure  of the
complexification  of $(\cS, I^{(z)})$, $z = U{\cdot}[0:1]$, along each horizontal leaf  $\cS \times \cS \times \{U\}$.
In other words, $\HCMqk$ is  the union of the complexifications of the  manifolds $(\cS, I^{(z)})$, $z \in S^2$.
This is  clearly not identical to a ``complexification" of $(\HMqk, \bI)$ in the sense of the previous paragraph;
rather it is merely a partial complexification in that context. 
Nevertheless, for the sake of brevity,  we  simply  call  it  the  {\it complexified harmonic space}. 
\par
\medskip
\section{Prepotentials for instantons on a qk manifold}
\setcounter{equation}{0}
\subsection{Lifts of  instantons}
Let  $(E, D)$  be a gauge field on a (strict)  qk manifold $(M, g, Q)$,   associated with the pair $(P, \o)$,  
a principal bundle $P$ over $M$ with a  connection $\o$. Further, the   {\it lift} $(P',   \o')$ of $(P, \o)$ to the 
harmonic space  $\pi^{\operatorname{qk}}\colon \HMqk {\to} M$  is the pair given by the lifted bundle  
$P'  = \pi^{{\qk}*} P$  over $\HMqk$ together with the  natural lift  $\o'$ to $P'$ of   $\o$   (see e.g. \cite[\S 2.4]{DPS}). 
The {\it lift of $(E,D)$}  to $ \HMqk$ is 
the gauge field $(E' =  \pi^{{\qk}*} E\,,\, D' =  \pi^{{\qk}*} D)$  given by 
\begin{itemize} [leftmargin=18pt]
\item[i)] the lifted vector bundle $q'\colon  E' =  \pi^{{\qk}*} E = P' \times_{G, \r} V \to \HMqk$ and
\item[ii)] the covariant derivative $D' =  \pi^{{\qk}*} D$ on $E'$ determined  by $\o' $.
\end{itemize}
 \par
 \begin{prop} \label{prop28} Let  $(E', D')$ be  a gauge field associated with a principal bundle with connection  $(P', \o')$ on  
 $\HMqk$.  Then 
   $(E', D')$ is the lift  of an instanton  on $M$ if and only if:
\begin{itemize} [itemsep=6pt plus 5pt minus 2pt, leftmargin=18pt]
\item[a)] The  $\o'$-horizontal vector fields  on $P'$ of  the infinitesimal transformations   of the  right $L$-action  
 on $\HMqk$  are  complete  and the group action generated by these vector fields 
has  no singular orbit. 
 \item[b)]   The curvature $F'$ of $D'$ is such that $F'(E, \cdot)  = 0$ for  any infinitesimal transformation $E$ of the 
   $L$-action. 
 \item[c)] For any $4n$-tuple  $(e'_{\pm a})$ of vectors  in the complex tangent spaces  $T^\bC_x \HMqk$, which  project onto  an adapted  complex  frame $(e_{+ a}, e_{-b})$ of $M$,  the following relations hold   
 \begin{gather}\label{3.10} F'(e'_{+a}, e'_{+ b}) = F'(e'_{-a}, e'_{- b})  = 0\ ,\quad F'(e'_{+a}, e'_{- b}) = - F'(e'_{-a}, e'_{+ b}) \ .
\end{gather}
 \end{itemize}
  \end{prop}
  \begin{pf}  By \cite[Prop.2.2]{DPS},  conditions (a) and (b) hold if and only if  $(E', D')$ is the lift of a gauge field $(E, D)$ on $M$.  
  Such a gauge field is an instanton if and only if the component $F^{(2)}$ of its curvature, defined in \eqref{2.7}, is identically zero. 
  This holds if and only if (c) holds.  \end{pf}
 \par
 \medskip
 \begin{cor} \label{cor5.2} The lift $(E^\cS, D^\cS)$ to the Swann bundle $\cS$ of an instanton $(E, D)$ on $(M, g, Q)$ is an instanton 
 on the hyperk\"ahler cone  $(\cS, h, \bJ_\a)$. \par
 If, in addition, the instanton $(E, D)$ is the complexification of an instanton with compact structure group $G^o$,  
 its lift $(E', D')$  to $\HMqk$ is real analytic and   admits  local holomorphic extensions on open sets of $\HCMqk$.
 \end{cor}
 \begin{pf} We recall that $\HMqk$ is  the harmonic space  of the hk manifold $(\cS, h, \bJ_\a)$, so the lift $(E', D')$ to $\HMqk$ 
  can also be  considered as the  lift  of a gauge field $(E^\cS, D^\cS)$  on the hk manifold  $\cS$. 
  The metric  $h$ and  the complex structures $\bJ_\a$ of $\cS$ imply the existence of an adapted complex frame 
  $(e^\cS_{+a}, e^\cS_{-b})$  for the  complexified  tangent spaces  $T^\bC_y \cS$, consisting of 
 \begin{itemize}
 \item[i)] four complex  vectors, say $e^\cS_{\pm 1}, e^\cS_{\pm 2}$,  vertical with respect to the projection  onto $M$; 
 in particular,  each of them  is a linear combination 
 of the projections from $\HMqk$ onto $\cS$ of  infinitesimal holomorphic  transformations of  $L$, and
 \item[ii)] additional $4n$ complex vectors $e^\cS_{\pm a}$, $ a \geq 3$,  projecting onto  the elements of an 
 adapted  complex  frame $(e_{\pm (a-2)})$  of  a  qk manifold $(M, t{\cdot} g, Q)$   homothetic to $(M, g, Q)$. 
 \end{itemize}
 Since  the curvature $F'$ of $(E', D')$ satisfies Proposition \ref{prop28} and   $(E', D')$  is   the lift of the gauge field $(E^\cS, D^\cS)$,    (i) and (ii) imply that $(E', D')$  satisfies all the hypotheses of \cite[Lemma 4.2]{DPS}  with respect to the base manifold $\cS$.  
 This implies that   $(E^\cS, D^\cS)$ is an instanton on $\cS$. 
 The final claim of the Corollary follows  from \cite[Prop. 4.1]{DPS} and    the fact that $(E', D')$ is  the lift of   $(E^\cS, D^\cS)$.  
 {\phantom{X}}\end{pf}\par
 \subsection{Adapted coordinates and adapted frames on harmonic space}
 \subsubsection{Holomorphically adapted  coordinates} \label{sectionadapted}
Let $x_o \in M$ be a point of a  qk manifold  and  $\cU$ one of  its neighbourhoods on  which there  exist trivialisations of 
the bundle $\pi\colon \cS \simeq  \sH^o {\to} M$ allowing  the  identifications $\cS|_{\cU} \simeq \cU \,{\times}\, \bH^*$.  
We may  write this as  $\cS|_{\cU} \simeq  \cU \times ( \bC^2 {\setminus} \{0\})$, representing  quaternions 
 $q = x^0 + x^1 {\bf i}  + x^2 {\bf j} + x^3 {\bf k}$ as  complex pairs,    
 $q = (\h= x^0 + {\bf i} x^1\,,\, \z = x^2 + {\bf i} x^3)$. In this way,  any   set of coordinates
 $(x^i)\colon \cU \subset M \to \bR^{4n}$ on $\cU$ induces corresponding  coordinates on $\cS|_{\cU}$,
 \beq \label{adaptedcoor} \wh \xi  = (x^1, \ldots, x^{4n};  \h, \z)\colon \sH^o|_{\cU} \longrightarrow  \bR^{4n}  \times (\bC^2 {\setminus} \{0\}),
 \quad (\h, \z) \neq (0,0) \ .\eeq
We call such a system of coordinates   {\it adapted coordinates} for $\cS$.
Now, denoting by $u^i_{\pm}$ the elements of the matrices    
$U = \begin{pmatrix} u^1_+ & u^1_-\\ u^2_+ & u^2_- \end{pmatrix}$ in $\SL_2(\bC)$, each element of the $L$-invariant  open set  
$\HMqk|_{\cU} = (\pi^{\text{qk}})^{-1}(\cU)\simeq  \cU  {\times}  (\bC^2 {\setminus} \{0\}) {\times}  \SL_2(\bC)$ 
is uniquely   determined  by a tuple of $4n$ real and $6$ complex numbers
 $$(x^1, \ldots, x^{4n}, \h, \z, u^i_{\pm})\qquad \text{constrained by}\qquad (\h, \z) \neq (0,0)\,,\  \det(u^i_{\pm}) =  1\,.$$
 They do not form a  true set of coordinates because they are  not functionally independent. 
 However, they are very convenient and, with a minor abuse of language,  we  call them the 
 {\it adapted coordinates}  of $\HMqk$.\par
 \medskip
 Finally,   consider the complex structure $I^{(z_o)}$, $z_o = [0:1] $,  of the horizontal leaf 
 $\cS|_{\cU}{\times} \{I_2\} \subset \HMqk$  (for the definition of  $I^{(z)}$, see  \S \ref{sect41}) and denote by  
 $(z^{j a})$, $j = 1,2$,   $1 \leq a \leq 2n+2$,   a set of coordinates on  the complexification $\cS|_\cU \times \cS|_{\cU}$ 
 of  the complex manifold  $(\cS|_{\cU} \simeq \cS|_{\cU}  {\times}  \{I_2\}, I^{(z_o)})$ satisfying the following conditions: 
 \begin{itemize}
 \item[1)] $z^{11} = \h$,  $z^{12} = \z$ and $z^{2 1} = \overline \h$, $z^{22} = \overline \z$.
 \item[2)] The   coordinates $z^{1a}$,  $ a\geq 3$,   and the   $z^{11}$  and $z^{12}$  restrict  to a set of $I^{(z_o)}$-holomorphic coordinates  on  $\imath(\cS|_{\cU}) \subset \cS|_U \times \cS|_{\cU}$. 
 \item[3)] The coordinates $z^{2a}$  with $a \geq 3$ restrict  to  the complex conjugates $z^{2a} = \overline{z^{1a}}$ on $\imath(\cS|_{\cU})$.  
 \end{itemize}
 For the remainder of $\HMqk|_{\cU} = \cS|_{\cU}  {\times}  \SL_2(\bC)$, i.e.\ the points with adapted coordinates  of the form  $(x^i, \h, \z, u^i_{\pm})$  with $U = (u^i_\pm) \neq  I_2$, we 
 choose  the $(4n + 8)$-tuple of  complex coordinates  
 $(z^{+a}, z^{-b}, u^i_{\pm})$ with
 $$z^{+a} \= u^+_j z^{j a}(x, \h, \z) \ ,\qquad z^{-a} \= u^-_j z^{j a}(x, \h, \z)\,. $$
Here  $u^\pm_j = \pm \e_{jk} u^k_\mp$  are  the entries of the inverse matrix  $(u^{\pm}_j) = (u^j_\pm)^{-1}$.
 Note that:
\begin{itemize}[itemsep=2pt plus 5pt minus 2pt, leftmargin=18pt]
\item[a)] For any given $U = (u^i_\pm) \in \SL_2(\bC)$, the restrictions of the  $z^{+ a}$ (resp. $z^{-a}$)  to  the  
horizontal leaf $\cS|_{\cU} {\times}  \{U\} \subset \HMqk$ are holomorphic (resp. anti-holomorphic) coordinates  with 
respect to the complex structure $I^{(z')}$ corresponding to the point $z' = U {\cdot} [0:1] \in \bC P^1 \simeq S^2$.  If the $z^{\pm a}$  are considered  to be independent, they are holomorphic coordinates for the complexification of the complex manifold   $(\cS|_{\cU}  {\times}  \{U\}, I^{(z)})$.
 \item[b)] At each point, the vector fields $\frac{\p}{\p z^{\pm 1}}$, $\frac{\p}{\p z^{\pm 2}}$ span the  complexified tangent space  of  the fibre   $\cS|_{x}$ of the Swann bundle over $x \in \cU \subset M$. 
 \end{itemize}
 We call the tuples  $(z^{\pm a}\!,\! u^i_{\pm})$  {\it holomorphically adapted  coordinates} of $\HCMqk$. \par
 \bigskip
\subsubsection{Infinitesimal $L$-transformations in holomorphically adapted coordinates}
Let $\gl = \Lie(L) =  \left( \bR + \su_2 \right) \ltimes \sl_2(\bC)$ and  denote by $\mathcal B$ the real basis for $\gl$,
  \beq \label{basisL} \mathcal B \= (G^o_D,  G^o_0, G^o_1, G^o_2,  H^o_0, H^o_{++}, H^o_{--}, i H^o_0, i H^o_{++}, i H^o_{--})\qquad \text{where}\eeq 
 \begin{itemize} 
 \item $G^o_D$  is the the  generator of  the abelian  subgroup $\{e^t, t \in \bR\}  \subset  \bH^*$; 
\item $(G^o_0, G^o_1, G^o_2)$ is the  basis of     $\su_2$,   
given by 
\beq 
G^o_0 \= \begin{pmatrix} i  & 0\\ 0 & - i \end{pmatrix}\ ,\qquad   G^o_1 \=  \begin{pmatrix} 0 & 1\\ -1 & 0 \end{pmatrix}\ ,\qquad G^o_2 \=  \begin{pmatrix} 0 & i\\ i & 0 \end{pmatrix}\ ; 
\eeq
\item  $(H^o_0, H^o_{++}, H^o_{--})$  is  the (complex) basis  of  $\sl_2(\bC)$, given by    
\beq\label{Hzero}  H^o_0 \= \begin{pmatrix} 1 & 0\\ 0 & - 1 \end{pmatrix}\ ,\qquad H^o_{++} \= \begin{pmatrix} 0 & 1\\ 0 & 0 \end{pmatrix}\ ,\quad H^o_{--} \= \begin{pmatrix} 0 & 0\\ 1 & 0 \end{pmatrix}\ .\eeq
\end{itemize}
Now, the isomorphism $\bH^*/\bZ_2\simeq \bR^* \times_{\bZ_2} (\Sp_1 \times_{\bZ_2} \Sp_{p,q})/ \Sp_{p,q}$ in 
\S \ref{Swannsection}  yields a Lie algebra isomorphism  
\beq \label{identification} \bH = \Lie(\bH^*/\bZ_2) \simeq  \bR + \su_2(\bC) = \gl/\sl_2(\bC)\ , \eeq
according to which  the unit quaternions $1$, ${\bf i}$, ${\bf j}$,  ${\bf k}$   of $\bH$ correspond to   (the equivalence classes  of)  $G^o_D$, $G^o_0$, $G^o_1$ and $G^o_2$, respectively.\par
\medskip

 For each element $G^o_\a$ or  $H^o_\b$  of  $\mathcal B$, we denote  by  $G_\a$ or  $H_\b$  the  corresponding  holomorphic vector field on $\HMqk$, given by  the   holomorphic part of   the infinitesimal transformation of the corresponding Lie algebra element. 
The  explicit expressions in adapted   coordinates of  the holomorphic vector fields $G_\a$ and $H_\b$ 
  can be   determined on the basis of  the diffeomorphism \eqref{tildej} and  the description of the corresponding right  actions,   given  in \S \ref{swann}.  We obtain
\beq \label{holvectors}
\begin{split}
 & H_{0}= u^i_+ \frac{\partial}{\partial u^i_+}  
                                   - u^i_- \frac{\partial}{\partial u^i_-} \ , \quad
H_{++}=  u^i_+ \frac{\partial}{\partial u^i_-} \ , \quad
H_{--}  =   u^i_- \frac{\partial}{\partial u^i_+}\ ,\\
 & G_D =  \h\frac{\p}{\p \h} + \z \frac{\p}{\p \z}\ , \qquad
    G_0 = i \h\frac{\p}{\p \h} + i \z \frac{\p}{\p \z} + u^i_+ \frac{\partial}{\partial u^i_+}  
                                   - u^i_- \frac{\partial}{\partial u^i_-} \ ,\\
 & G_1 = - \z \frac{\p}{\p \h} +     u^i_+ \frac{\partial}{\partial u^i_-} -    u^i_- \frac{\partial}{\partial u^i_+}\ , \quad
  G_2 = i \z \frac{\p}{\p \h}  +     i u^i_+ \frac{\partial}{\partial u^i_-} +  i  u^i_- \frac{\partial}{\partial u^i_+}\,.
\end{split}    
\eeq
The  Lie brackets amongst these vector fields  are then easily computable:
\begin{align}
\nonumber & [G_0, G_1] = G_2\ , && [G_0, G_2] = - G_1\ ,&& [G_1, G_2] = G_0\ ,\\
\nonumber & [H_0, H_{++}] = 2 H_{++}\ , && [H_0, H_{--}] = - 2 H_{--} \ ,\!\!&& [H_{++}, H_{--}] = H_0\ ,\\
\nonumber &[G_0, H_0] = 0\ ,&& [G_0, H_{++}] =  2 i H_{+ +}\ ,&& [G_0, H_{--}] = \pm 2 i H_{--}\ ,\\
\nonumber  &[G_1, H_0] = - 2 (H_{++} {+} H_{--}), && [G_1, H_{++}] =  H_0\ ,&& [G_1, H_{--}] =  H_0\ ,\\
 \nonumber&[G_2, H_0] = - 2 i (H_{++} {-} H_{--} ),\!\!\!&& [G_2, H_{++}] =  - i H_0\ ,&& [G_2, H_{--}] =  i H_0\ , \\
& [G_D, G_\a] = 0\;, \ \ \a = 0,1,2\ ,&&  [G_D, H_\b] = 0\;,&& \hskip -1.3cm  \b \in \{ 0, \pm \pm\}\ .
\end{align}

\par
\medskip
\subsubsection{Strongly adapted  frames of  $(\cS, h, \bJ_\a)$}
In the above collection of holomorphic vector fields, the $H_\b$  are vertical with respect to $\pi^{\cH|\cS}: \HMqk \to \cS$, 
whereas  the $G_A$  uniquely decompose as  sums  $G_A = G_A^\cS + V_A$,  with $V_A$  vertical and  $G^\cS_A$  
tangential to the horizontal leaves $\cS \times \{U\} \subset \HMqk $.  This decomposition of the  $G_A$  can be 
explicitly determined by examining  \eqref{holvectors}.  
There, we can also see that, denoting by  $\k$ the constant   $\k = |\L|$  occurring in the definition   \eqref{swannmetric}
of the hyperk\"ahler metric $h$,  the  vector fields
 \beq 
\begin{split} &f_{11} \= \frac{1}{2 \k(|\h|^2 + |\z|^2)^{\frac{1}{2}}} \left(G_D^\cS - i G^\cS_0\right)=  \frac{1}{\k(|\h|^2 + |\z|^2)^{\frac{1}{2}}} \left(\h\frac{\p}{\p \h} + \z\frac{\p}{\p \z}\right), \\
 & f_{12} \= \frac{1}{2\k(|\h|^2 + |\z|^2)^{\frac{1}{2}}} \left( G_1^\cS - i G^\cS_2\right) =   \frac{1}{\k(|\h|^2 + |\z|^2)^{\frac{1}{2}}}\left(- \z \frac{\p}{\p \h} + \h \frac{\p}{\p \z}\right)
\end{split}
 \eeq
 give   $h$-orthonormal   frames  at each vertical  tangent space fibre $T^V \cS_x$, $x \in M$, of the Swann bundle. 
 \par
  We  recall that  the hyperk\"ahler structure $(\bJ_1, \bJ_2, \bJ_3)$ of $\cS$ acts on vectors as  the infinitesimal transformations of  $\bH = \Lie(\bH^*/\bZ_2) $  associated with   the elements ${\bf i}, {\bf j}, {\bf k} \in \bH$.  Thus, by  the identification  \eqref{identification},  we have  
  $$\bJ_1 G^S_D = G^\cS_0\ ,\qquad \bJ_2 G^S_D = G^\cS_1\ ,\qquad  \bJ_3 G_D = G^\cS_2 $$
and now setting   $f_{2j} \= \overline{f_{1j}}$, $j = 1,2$,  we see that 
\beq \bJ_1 f_{1j} =  i f_{1j}\ , \quad \bJ_1 f_{2j} = -  i f_{2j} \ ,\quad  \bJ_2 f_{11} =   f_{22} 
  \ ,\qquad  \bJ_2 f_{12} = -  f_{21}\ . \eeq
Next we want to complete this set of four vector fields $f_{ij}$, $i ,j = 1, 2$,  to a {\it strongly adapted} complex frame field. In order to do this, we consider the horizontal leaf $\cS|_{\cU} \times \{I_2\} \subset \HMqk$  and  a field of adapted complex frames
$(e_{\pm a})$ for  $T^\bC \cS \simeq T^\bC ( \cS \times \{I_2\})$,  in which the elements $e_{\pm a}$ with $a = 1,2$   are given by precisely the fields $f_{ia}$,
\beq e_{+ 1} = f_{1 1} \ ,\quad e_{+2} = f_{1 2}\ , \quad  e_{-1} = f_{2 1} \ ,\quad e_{-2} = f_{2 2}\ .\eeq
Then, over  any other  horizontal leaf,  $ \cS|_{\cU}  \times \{(u^i_{\pm})\} \subset  \HMqk$,  the adapted complex frame field $(e_{\pm a})$  is given by
$e_{\pm a} = u^j_\pm e_{i a}$. In particular, by construction,
     \beq \label{stronglyad} e_{+ 1} = u^i_+ f_{i 1} \ ,\quad e_{+2} = u^i_+ f_{i 2}\ , \quad  e_{-1} = u^i_- f_{i 1} \ ,\quad e_{-2} = u^i_- f_{i2}\,.\eeq
\par
\subsection{Prepotentials and  instantons}
\subsubsection{The   components of  the potential of  a  lifted instanton}
Let $(E', D')$ be the lift  to $\HMqk$ of an instanton $(E, D)$ on a qk manifold $(M, g, Q)$.  
Consider  an open subset   $\cW$ of $\HMqk$ with  holomorphically adapted   coordinates $(z^{\pm a},  
u^i_{\pm})$  and  gauges (trivialisations)  for the lifted  $G$-bundle $P'$ and the associated vector bundle $E'$, 
 $$ \f\colon P'|_{\cU} \to \cW \times G\ ,\qquad \wh \f\colon  E'|_{\cU} \to \cW \times V\ .$$
   By Corollary \ref{cor5.2},   $(E',  D')$ is the lift of an  instanton $(E^\cS, D^\cS)$ on the hk manifold $\cS$. Then, by  \cite[Prop. 4.1]{DPS},  the bundles $P'$ and $E'$ are real analytic, they  have  natural structures of holomorphic bundles over $\HMqk$ and the connection $1$-form $\o'$ (corresponding to $D'$) is invariant under the complex structure of $P'$.  We may therefore assume  that the gauge $\f$ is  holomorphic and that in this gauge $\o'$ is  determined by its $(1,0)$-potential $A^{10}$, the restriction of the potential $A$ to  the holomorphic tangent space $T^{10} \HMqk$ of $\HMqk$.  
Because of real analyticity, we may also assume that the components of $A$ admit  unique holomorphic extensions  from  $\cW$ to a    neighbourhood  $\wt\cW$ of  complexified harmonic space   $\HCMqk = \HMqk \times \HMqk$. 
      \par
      \medskip
We now recall that  at  each fixed  $ (z^{\pm a}_o, u^i_{o\pm}) \in \wt \cW$,   the vectors   
$$ H_0|_{y_o}\ ,\qquad H_{\pm\pm}|_{y_o}\,,\qquad \frac{\p}{\p z^{\pm a}}\bigg|_{z^{\pm a}_o}$$
constitute a basis  for the  holomorphic tangent space $T^{10}_{y_o} \HCMqk$. 
 This means that the $(1,0)$-potential $A^{10}$ of  $(E', D')$ is uniquely determined by the $\gg$-valued holomorphic functions 
  \beq \label{nota} A_0 \= A(H_0)\ ,\qquad  A_{\pm\pm} \= A(H_{\pm\pm})\ ,\qquad  A_{\pm a}\= A\left(\frac{\p}{\p z^{\pm a}}\right)\ .
  \eeq
 In  the   following, $A_0$, $A_{\pm\pm}$,  $A_{\pm a}$  will  always denote these  $\gg$-valued functions. \par 
\par
\subsubsection{The analytic gauge and the prepotential of  an  instanton}
Keeping to the  notation of the previous subsection, pick a point    $ (z_o^{\pm a},  U_o) \in \cW$ with 
$U_o = I_2$  and write  $x_o = \pi^{\text{qk}}(  z_o^{\pm a},  U_o) ) \in M$. 
The following theorem directly follows from  \cite[Thm.4.3 and Lemma 4.6]{DPS}.\par
\begin{theo} \label{thm53}  
The neighbourhood $\wt \cW \subset \HCMqk$  and the   holomorphic gauge 
$\varphi\colon P'|_{\wt \cW} \to \wt \cW \times G$ can be chosen in such a way  that  the components  of the instanton potential $A^{10}$ satisfy the conditions 
\beq \label{claim1} A_0 =  A_{-a}  = A_{+1} = A_{+2} = 0 \ .
\eeq
Moreover,
\begin{itemize}[itemsep=2pt plus 5pt minus 2pt, leftmargin=18pt]
\item[1)] the component $A_{--}$  satisfies the differential conditions 
\beq\label{5.15bis}  H_0{\cdot} A_{--} = - 2 A_{--}\ ,\qquad \frac{\p A_{--}}{\p z^{-a}} = \frac{\p A_{--}}{\p z^{+1}} = \frac{\p A_{--}}{\p z^{+2}}  = 0 \eeq
\item[2)] the   components $A_{++}$,  $A_{+a}$  are   uniquely    determined  by   $A_{--} $, being the unique solutions of  the differential problem  
\beq   \label{5.15}
\begin{split}
&  H_0{\cdot} A_{++} =  2 A_{++}\ ,  \qquad H_{--}{\cdot} A_{++} - H_{++}{\cdot} A_{--}  -[A_{++}, A_{--}] = 0\ ,\\
& A_{+a} = - \frac{\p A_{++} }{\p z^{-a}} \ . 
\end{split}
\eeq
\end{itemize}
\end{theo}
\begin{pf}  The proof  of the  existence of  a holomorphic gauge $\varphi\colon P'|_{\wt \cW} \to \wt \cW \times G$ 
in  which $A^{10}$ satisfies \eqref{claim1} involves a slight modification to  the proof of  \cite[Thm. 4.3] {DPS}.  Consider 
the $G$-invariant distributions $\cD_-^h$, generated by the horizontal lifts of the vector fields $\frac{\p}{\p z^{-a}}$, 
and $\cD_-^h \oplus \langle H^h_0, \left(\frac{\p}{\p z^{+1}}\right)^h, \left(\frac{\p}{\p z^{+2}}\right)^h\rangle$, 
where the second summand is generated by  the horizontal lifts of $H_0$, $\frac{\p}{\p z^{+1}}$ and $\frac{\p}{\p z^{+2}}$.
Both these distributions are  involutive, the second because
\begin{itemize}[itemsep=2pt plus 5pt minus 2pt, leftmargin=18pt]
\item[a)]  $P'$, being the lift of an instanton on $M$, has identically vanishing curvature along the infinitesimal transformations of the structure group $L$ of $\HMqk$ over $M$ (see e.g. \cite[Prop.2.2]{DPS}), and
\item[b)] both  $\frac{\p}{\p z^{+1}}$ and  $\frac{\p}{\p z^{+2}}$ are pointwise linear combinations of the infinitesimal transformations $H_\a$ and $G_\b$. 
\end{itemize}
We may now follow the argumentation in the proofs of \cite[Thm 4.3 and Lemma 4.6]{DPS} to obtain all the
claimed results except for $ \frac{\p A_{--}}{\p z^{+j}}  = 0\,,\,j = 1,2$. These equations, however, follow from other   
conditions in  \eqref{5.15} and  \eqref{claim1}:
\beqn
\begin{split}  0 &= \frac{\p}{\p z^{-j}}{\cdot} \left(H_{--}{\cdot} A_{++} - H_{++}{\cdot} A_{--}  -[A_{++}, A_{--}] \right) \\
&= H_{--}{\cdot} \left( \frac{\p A_{++}}{\p z^{-j}} \right) +  \frac{\p A_{--}}{\p z^{+j}}- \left[ \frac{\p A_{++}}{\p z^{-j}}, A_{--}\right] \\
&= - H_{--}{\cdot} A_{+j} +  \frac{\p A_{--}}{\p z^{+j}}+  \left[ A_{+j}, A_{--}\right]   =  \frac{\p A_{--}}{\p z^{+j}}\ .
\end{split}\eeqn
\vskip-0,7cm \end{pf}

A holomorphic gauge $\f\colon P' \to \wt \cW \times G$ in which the  potential $A^{10}$ satisfies \eqref{claim1} 
is called an {\it analytic gauge}.   The component $A_{--}$ of  $A^{10}$ in  an analytic  gauge is  the {\it prepotential} 
of the instanton in that gauge .\par
\medskip
The above theorem essentially says that an instanton $(E, D)$ on a qk manifold is completely determined its prepotential
$A_{--}$  in an analytic gauge  for  the lifted gauge field $(E', D')$.  For  hk manifolds  this property admits
a  converse  \cite[Thm. 5.1]{DPS}.  We now give the corresponding property for  qk manifolds.

\begin{theo} \label{thm54}
Let $\cW\subset  \HMqk$  be an open subset of the harmonic space of a qk manifold $(M, g, Q)$ on which we may 
choose  holomorphically adapted  coordinates $(z^{\pm a}, u^i_{\pm})\ ,\ a=1,\dots,2n+2\,,\,i=1,2$ and  a trivialisation 
of the lifted bundle $P'$ over  $\HMqk$. Given a  map  $\bA_{--}\colon \cW \subset \HMqk \to \gg$ satisfying  
\beq \label{prep-1} 
H_0{\cdot} \bA_{--} = - 2 \bA_{--}\ ,\qquad \frac{\p \bA_{--}}{\p z^{-a}} = 0\ ,\qquad  
\frac{\p \bA_{--}}{\p z^{+1} }= \frac{\p \bA_{--}}{\p z^{+2}} = 0
\eeq
there exists a unique (up to equivalence) instanton $(E, D)$ on the qk manifold 
$(\cU\=  \pi^{\operatorname{qk}}(\cW) \subset M$, $g|_{\cU}, Q|_{\cU})$,  admitting an analytic gauge  in which 
$\bA_{--}$ is the prepotential. 
\end{theo}
\begin{pf} By \cite[Thm. 5.1]{DPS} and the first two equations in \eqref{prep-1} we have existence of 
an instanton $(\wh E, \wh D)$ on the open subset $\cU^\cS = \pi^{\cH|\cS}(\wt \cW)$ of the  hyperk\"hler cone 
$(\cS, h, \bJ_\a)$ with  prepotential  $A_{--} = \bA_{--}$.  The claim is then proven if  $(\wh E, \wh D)$ is shown
to be the lift of a gauge field on  $\cU = \pi^{\operatorname{qk}}(\cW)$ of $M$. By \cite[Prop.2.2]{DPS}, this is the
case if and only if two conditions are satisfied: a) the horizontal lifts  to $\wh P$ of the  infinitesimal transformations 
$G^\cS_\a$, $\a \in \{D, 0, 1, 2\}$,  of the structure group $\bH^*/\bZ_2 = L/\SL_2(\bC)$  of the Swann bundle over $M$ 
are complete vector fields and  b) the curvature $F$ of  $(\wh E, \wh D)$ vanishes identically when evaluated along 
any $G^\cS_\a$. Condition a) is surely  satisfied: For any trivialisation $P'|_{\cU} \simeq \cU \times G$,  the $\cU$-component of the horizontal lift  of the vector field $G^\cS_\a$  is  $G^\cS_\a$ itself and the  $G$-part  
is the image of  $G^\cS_\a$  in $\gg$  by applying to it the potential of the connection in this trivialisation.
Both components are therefore  bounded at all points and the horizontal lifts of the $G^\cS_\a$ are therefore complete. 
It remains to check that  b) holds. 
Since the vector fields $G^\cS_\a$   are pointwise linear combinations of the complex vector fields 
$\frac{\p}{\p z^{\pm j}}$, $j = 1,2$, we only need to check that 
$$F\left(\frac{\p}{\p z^{+ j}}, \frac{\p}{\p z^{- a}}\right) = F\left(\frac{\p}{\p z^{- j}}, \frac{\p}{\p z^{+ a}}\right)  = 0\,,\quad  
 j = 1,2\ ,\  a = 1, \ldots, 2n +2\ . $$
By  \cite[(5.19)]{DPS}, this holds if and only if 
\beq \label{tempo1} \frac{\p^2 A_{++}}{\p z^{- j} \p z^{-a}}  = 0\ .\eeq
Now, recall (see \cite[(5.20)]{DPS})  that  the component $A_{++}$  of the lift of $(\wh E, \wh D)$ to $\cW$ satisfies the differential equation
$$H_{--}{\cdot} A_{++} - H_{++}{\cdot} \bA_{--} + [\bA_{--}, A_{++}]  = 0\ .$$
Differentiating this equation along $\frac{\p}{\p z^{-j}}\,,\,j=1,2$, and  using \eqref{prep-1}, we obtain
\beq \label{tempo} 
 H_{--}{\cdot}  \frac{\p A_{++}}{\p z^{-j}}     +  \left[ \bA_{--}\, \frac{\p A_{++}}{\p z^{-j}} \right] = 0\ .\eeq
Thus, the $G$-equivariant extension of the functions $g_{+j} \= \frac{\p A_{++}}{\p z^{-j}}$ on $P'$ satisfy the  
differential problem 
 $$H^h_0{\cdot} g_{+j} = g_{+j} \ ,\qquad  H^h_{--}{\cdot} g_{+j} = 0\,,$$
with  $H^h_0$ and  $H^h_{--}$ being horizontal lifts to  $P'$ of $H_0$ and  $H_{--}$, respectively.  
Now, using \cite[Lemma 5.3]{DS}, in the version obtained by replacing  $H_{++}$  by  $H_{--}$, shows that 
this system only has a trivial solution, $\frac{\p A_{++}}{\p z^{- i} }  = g_{+i}|_{\HMqk \times \{e\}} = 0$,
which yields the required condition  \eqref{tempo1}.
\end{pf}
 
\subsubsection{Equivalent  prepotentials and curvatures} \label{subsub53}
Let  $\cW\subset  \HMqk$  be an open subset of the harmonic space of   $(M, g, Q)$ admitting a set of  holomorphically 
adapted  coordinates $(z^{\pm a}, u^i_{\pm})$ and  a trivialisation of the lifted bundle $P'$ associated with a vector bundle  $E$.  \par
As shown by  Theorems \ref{thm53} and \ref{thm54},  instantons $(E, D)$ are completely  determined by their  
prepotentials, the real analytic maps $A_{--}\colon \cW \to \gg$ satisfying \eqref{prep-1}. Two distinct $\gg$-valued 
prepotentials $A_{--}$, $A'_{--}$  on $\cW$ correspond two equivalent instantons if and only if they are the 
$H_{--}$-components of two $(1,0)$-potentials $A^{10}$, $A'{}^{10}$, both in an analytic gauge (i.e. satisfying 
\eqref{claim1}) and related to each other by a gauge transformation  \eqref{changepot} with $g\colon \cW \to \gg$.  

On the other hand, if the exponential map of the structure group $G$ is surjective (as e.g. when $G$ is reductive),  any gauge transformation $g\colon \cW \to G$ can be written in the form $$g(z^i, u^j_{\pm}) = e^{\wh g(z^i, u^j_{\pm})}\qquad \text{for some appropriate}\ \wh g\colon \cU \to \gg\ .$$
In this case,  {\it $A_{--}$,  $A_{--}'$  correspond to  equivalent instantons if and only if 
\beq \label{equiv}
\begin{split}
&A'_{--} = \Ad_{e^{-\wh g}}(A_{--} + H_{--}{\cdot} \wh g)\\
& \text{for some}\  \wh g\colon \cU \to \gg\ \ \text{such that}\qquad  H_0{\cdot} \wh g =  \frac{\p \wh g}{\p z^{-a}} =  \frac{\p \wh g}{\p z^{+1}}  =  \frac{\p \wh g}{\p z^{+2}}  = 0
\end{split}
\eeq
In particular,  $A_{--}$ corresponds to a flat connection if and only if $A_{--} = - H_{--}{\cdot} \wh g$ for some  $\gg$-valued map $\wh g\colon \cW \to \gg$ satisfying \eqref{equiv}}.\par
\par
\medskip
We now  recall  that  the curvature of the  $(1, 0)$ potential $A^{10}$ of an instanton on $\cW\subset  \HMqk$ in an 
analytic gauge is determined by the  $A_{++} = A^{10}(H_{++})$  component by the simple formula \cite[(5.19)]{DPS},
\beq\label{5.7bis} \cF^\f\left(\frac{\p}{\p z^{+a}}, \frac{\p}{\p z^{-b}}\right) = \frac{\p^2 A_{++}}{\p z^{+a} \p z^{+b}}\ .\eeq 
Due to this simple relation with the curvature,  the component $A_{++}$ may be thought of as a {\it second prepotential} of the instanton $(E, D)$ on $M$ \cite{DPS}. By Theorem \ref{thm53}, the  potential $A^{10}$, and thus  the component $A_{++}$,  is uniquely determined by the prepotential $A_{--}$.  A direct  determination of $A_{++}$ from $A_{--}$
follows from the following  (cf.  \cite[Prop. 5.5]{DPS}):
\begin{prop} \label{prelim} Given a prepotential  $A_{--}\colon  \cW \subset \HMqk  \to \gg$ for a lifted instanton in an analytic gauge, the corresponding second prepotential
$A_{++}$ coincides with  the  {\rm unique}   solution on $\cW$ of the differential problem for
the unknown  $B_{++}$,
\beq \label{characterising}  H_{--}{\cdot} B_{++} = H_{++}{\cdot} A_{--} -[A_{--}, B_{++}]\ ,\qquad H_0{\cdot} B_{++} = 2 B_{++}\ .
\eeq
\end{prop}
\noindent The horizontal projection of $B_{++}$  in the bundle  $\HMqk\to M$  is known as the {\it Leznov field} and plays a r\^ole in various contexts  (see e.g.\ \cite{Lez,DL2, Si,DL1,DO,De}).
\par
\medskip
\section{Prepotentials in canonical forms} 
\subsection{Potentials in canonical form  for  gauge fields of compact type}
\subsubsection{The canonical connection of a Stiefel bundle}
 Let  $(e_1^o, \ldots, e_{m+k}^o)$ \label{Stiefel1}  be the  standard  basis of   $\bC^{m+k}$ and    denote by   $\cV_{m+k,m}$ the  {\it Stiefel manifold} of  unitary $m$-frames  in  $\bC^{m+k}$, i.e. the manifold of all  ordered $m$-tuples $(v_1, \ldots, v_m)$ of complex vectors 
 $v_i = b_i^\a e_\a^o$ given  by  some  $B \= (b_i^\a) \in \bC_{(m+k) \times m}$  such that
\beq \label{condA}  B^\dagger  B = \Iden_{m}\ ,\qquad  	\text{where}\ \  B^\dagger := \overline{B^T}\,.\eeq
Note that for any  such unitary $m$-frame,  the  associated linear map
\beq P_{B}\colon  \bC^{m+k} \to \bC^{m+k}\ ,\qquad  P_{B}\left(\begin{array}{c}x^1\\ \vdots\\ x^{m+k}\end{array}\right) \= B  B^\dagger \left(\begin{array}{c}x^1\\ \vdots\\ x^{m+k}\end{array}\right)\eeq
is the unitary projector  from $\bC^{{m+k}}$ onto  the $m$-plane spanned by the  $v_i$. In particular,  $P^2_B = P_B$.
The manifold $\cV_{{{m+k}},m}$ is    identifiable with the homogeneous space
$\cV_{{{m+k}},m} = \U_{{m+k}}/(\{\Iden_m\} \times \U_{k})$  and its natural projection onto  the Grassmannian of $m$-planes 
\beq \label{Stiefelbundle} 
\pi\colon \cV_{m+k, m} \to \operatorname{Gr}_m(\bC^{{m+k}})  \simeq  \U_{{m+k}}/(\U_m \times \U_{k})\eeq
makes it  a   principal $\U_k$-bundle over  $\operatorname{Gr}_m(\bC^{{m+k}})$. 
Consider the map  $\cM\colon \cV_{{{m+k}},k}\to  \bC_{{{(m+k)}} \times k}$ which sends  each $m$-frame 
$$(v_1 \= B_1^\a e^o_\a, \ldots, v_m \= B_m^\a e_\a^o)$$
 into the corresponding  matrix $B = (B_i^\a)$.  The    $\bC_{m\times m}$-valued $1$-form  on $\cV_{{{m+k}},m}$ 
\beq \o^{\text{can}} =  \cM^\dagger  d \cM\eeq
  is then  $\gu_m$-valued.   Further, it   can be checked that it is   a connection $1$-form on the  $\U_m$-bundle    \eqref{Stiefelbundle}, called   the {\it canonical connection} of $\cV_{{{m+k}},m}$ \cite{NR}. 
\par
\medskip
\subsubsection{Potentials in canonical form}
Let  $G^o$  be  a compact Lie group, which,  with no loss of generality, we  assume to be  a closed subgroup of a unitary group $\U_m$  for some integer $m$. The following crucial fact has  been proved by Narasimhan and Ramanan.\par
\begin{theo}[\cite{NR}]  \label{NR}  Let  $(E^o, D^o)$   be a  gauge field of compact type  over  an $n$-dimensional  manifold $M$,  associated with a principal bundle with connection $(P^o, \o^o)$ having structure group $G^o \subset \U_m$. Then  there exists a $G^o$ equivariant embedding 
\beq \label{TheoNR} \Phi\colon P^o \to \cV_{{m+k}, m}\qquad \text{with}\qquad k= (n+1)(2n+1) m^2 - m\ ,\eeq
    such that 
$\o^o = \Phi^*\left(\o^{\text{\rm can}}\right) =   \cM^{\Phi\dagger}d \cM^{\Phi} $ with  $\cM^{\Phi} = \cM \circ \Phi$.
\end{theo}
This  implies  that,  in a fixed   trivialisation  $P^o|_{\cU} \simeq \cU \times G^o$, the  potential $A$  of the connection $1$-form $\o^o$ can always be written as 
\beq\label{universal} A = \x^\dagger d \x\ ,\qquad \text{where}\qquad \x\= \left(\cM^{\Phi}\right)|_{\cU \times \{e\}} = \colon \cU \to \bC_{(m+k) \times m}\ .\eeq
Note that  by construction,    $\x^\dagger_x \x_x = \Iden_m$ at  any $x \in \cU$. Observe  also that: \\[10pt]
 (1) If $A'$ is the potential of $\o^o$ in a new gauge,  related to $A$  by  a gauge transformation determined by 
 a map  $g\colon \cU \to G^o \subset \U_m$, then   $A'$  also has  the form  \eqref{universal}, with  $\xi$ is replaced by  
 $ \xi' = \xi g$. This is due to the fact that 
\beq \label{Uform} A' = \Ad_{g^{-1}} A  + g^{-1}  dg =   (g^\dagger \xi^\dagger) (d \xi g)   +  (g^\dagger \xi^\dagger) (\xi  dg) =  (\xi g)^\dagger d  ( \xi g)\ .\eeq
 \noindent
(2) From the definition of the map  $\xi =\cM^\Phi|_{\cU \times \{e\}}$, after possibly restricting the open set $\cU$, 
it is always possible to determine  a
 {\it row-permuting matrix}  $C \in \U_{m+k}$  such that  the upper $m \times m$  block
$h_x$ of the product    $C  \xi_x \in \bC_{(m+k) \times m}$ is  of maximal rank for all $x \in \cU$.   
This determines a map  $h\colon \cU \to \GL_m(\bC)$   such that  $C \xi h^{-1}$ has the form 
\beq \label{lambda} C  \x  h^{-1} = \left(\begin{array}{c} \Iden_m \\  \l \end{array}\right)\qquad \text{for some map }  \l\colon \cU \to \bC_{k \times m} \eeq
 uniquely determined  by  $\xi$ and the constant matrix $C$. \\[5pt]
(3) Conversely, given a smooth map $\l\colon \cU \to \bC_{k \times m}$ we can find a row-permuting matrix $C$ and a map
$\xi\colon \cU \to \bC_{(m+k) \times m}$  which determine $\l$. For example, $C = \Iden_{m+k}$ and 
\beq \label{canonicalxi}  \xi = \xi^{(\l)} \= \left(\begin{array}{c} \Iden_m \\  \l \end{array} \right)  (\Iden_m + \l^\dagger \l)^{-\frac{1}{2}}\eeq
will do. \\[5pt]
(4) Consider two smooth maps $\l, \l'\colon \cU \to \bC_{k \times m}$ and the  $\U_k$-valued maps 
$$A^{(\l)} = \xi^{(\l)\dagger} d \xi^{(\l)} \ ,\qquad A^{(\l')}= \xi^{(\l')\dagger} d \xi^{(\l')}$$
with $\xi^{(\l)}$, $\xi^{(\l')}$  as in \eqref{canonicalxi}. 
The maps $A^{(\l)} $ and $A^{(\l')} $  are two  potentials   of the same connection   in two   distinct gauges if and only if 
there exist two  row-permuting  matrices $C$, $C' \in \U_{m+k}$ and  a  map  $g\colon \cU \to G^o \subset \U_m$,  
such that  $C' \xi^{(\l')}= C \xi^{(\l)} g$.   On the other hand,  since  $C^\dagger = C^{-1}$ and  $C'{}^\dagger = C'{}^{-1}$, 
we see that  $A^{(\l)}$ and $A^{(\l')}$ are independent of  $C$ and $C'$ and, consequently, that there is  no loss of 
generality  in assuming  that $C = C' = \Iden_{m+k}$.   Hence, $A^{(\l)} $ and $A^{(\l')} $ are gauge equivalent 
if and only if 
$$ \left(\begin{array}{c} \Iden_m \\  \l \end{array} \right)  = \left(\begin{array}{c} \Iden_m \\  \l' \end{array} \right) \wt g \qquad \text{with}\qquad  \wt g =  (\Iden_m + \l^\dagger \l)^{\frac{1}{2}} (\Iden_m + \l'{}^\dagger \l')^{-\frac{1}{2}}  g\ .  $$
This holds if and only if $\ \wt g = \Iden_m$ and $\l = \l'$, from which $g = \Iden_m$ follows.\par
\smallskip
These observations imply  that  in each  gauge equivalence class of  potentials, there exists  {\it  exactly one potential of the form $A^{(\l)} = \xi^{(\l)\dagger} d \xi^{(\l)} $ for some map $\l\colon \cU \to \bC_{(m+k) \times m}$}. 
The number of distinct $\l$'s which determine the same potential  $A^{(\l)}$ is less than or equal to $(m+k)!$,  the cardinality of the group of  row-permuting  matrices $C$. \par
\medskip
We have therefore verified  the following direct  consequence of Theorem \ref{NR}. 
\begin{cor} \label{lemma62} For any  gauge field  $(E^o, D^o)$   with a structure group $G^o \subset \U_m$,  there exists an integer $k \geq 0$ with the property  that, around   any  $x_o \in M$,  there is a   neighbourhood $\cU$ and a  gauge   $\f\colon P^o|_{\cU} \to \cU \times G$,  in which   the  potential  $A\colon T^* \cU \to \gg^o $ has the form  
\beq \label{canonicalform}  A = A^{(\l)} = \xi^{(\l)\dagger} d \xi^{(\l)}\eeq
with $\xi^{(\l)}$    as in \eqref{canonicalxi} for  some smooth $\l\colon \cU \to \bC_{k \times m}$.  Such a gauge  is uniquely 
determined  by the restricted gauge field $(E^o, D^o)|_{\cU}$ and  the map  $\l$ is uniquely determined 
up to the  action of a  row-permuting matrix $C$,  according to the equations \eqref{lambda}  and \eqref{canonicalxi}.
\end{cor}
\noindent
The potential \eqref{canonicalform}  is  said to be in    {\it canonical form} and we call the corresponding map  
$\l\colon \cU \to \bC_{k \times m}$  its  {\it core}.
\medskip
\par
\subsection{The expression of   the curvature in terms of the core }
Let $A^{(\l)} = \xi^{(\l)\dagger} d \xi^{(\l)}\colon \cU \to \gg^o$ be a potential in canonical form with core $\l$  for a gauge field $(P^o, \o^o)$ with structure group $G^o$. 
For  such potential, we  now  select a smooth map  $\h\colon \cU \to \bC_{(m+k)\times k}$ with the following  property: 
{\it for each $x \in \cU$   the set of columns of the matrix $\h_x$ complete the set of columns of  $\xi^{(\l)}_x$ to a unitary basis for the vector space $\bC_{(m+k)\times (m+k)}$}.  In other words, we assume  that  $\h$ is such that 
\beq
\text{the matrix}\ \ \left(\xi^{(\l)}_x\big| \h_x\right) \ \ \text{is in}\ \ \U_{m+k}\qquad \text{for each}\ x \in \cU\ .\eeq
From the above characterisation, we clearly have that  $\h^\dagger \h = \Iden_k$. Moreover, since  the columns of $\h_x$ are linearly independent from those of $\xi_x$,  we may always  assume that $\h$ has the  form
\beq
 \label{eta}
 \h =  \left(\begin{array}{c} \nu \\  \Iden_k \end{array}\right) (\Iden_k + \nu^\dagger \nu)^{-\frac{1}{2}}\qquad \text{where} \ \nu\colon \cU \to \bC_{m \times k}\,. 
\eeq

Note that the columns of such a map $\h$ are orthogonal to all the columns of $\xi^{(\l)}$ if and only if
\beq \nu^\dagger + \l = 0\ , \eeq
a condition which uniquely determines $\h$ in terms of  $\l$.  
In what follows, we denote  the map  \eqref{eta}  by $\h^{(\nu)}$ and we say that  $\nu = - \l^\dagger$ is the  {\it conjugate core} of the potential $A = A^{(\l)}$.\par
Clearly {\it a  core $\l$ and the associated map $\xi^{(\l)}$  are uniquely determined by the conjugate core $\nu = - \l^\dagger$ and the corresponding map $\h^{(\nu)}$, and vice versa.}\par
The following lemma (see e.g.  \cite[Ch. II.3]{At}) gives a simple and convenient expression for the curvature of a connection in terms of the two  maps $\xi^{(\l)}$ and $\h^{(\nu)}$.
\begin{lem} \label{curvature1}  Let   $A^{(\l)} = \xi^{(\l)\dagger} d \xi^{(\l)}\colon \cU \to \gg^o$ be a potential  in canonical form with  core $\l$ and $\nu = - \l^\dagger$.  The curvature 
$F$ of the  corresponding  gauge field  is
\beq \label{curvature11}  F =   \xi^{(\l)\dagger}  \left\{ \left(\begin{array}{cc} d \nu\wedge   (\Iden_{2k} + \nu^\dagger \nu)^{-1 } d \nu^\dagger & 0\\ 0 & 0  \end{array} \right)\right\} \xi^{(\l)}\,. \eeq
\end{lem}
\begin{pf}  Since $\xi^{(\l) \dagger} \xi^{(\l)} = \Iden_m$ and $\xi^{(\l) \dagger} d \xi^{(\l)} = - d \xi^{(\l) \dagger}\xi^{(\l)}$,  
\begin{align} \label{615}  \xi^{(\l)\dagger} \left\{d (\xi^{(\l)} \xi^{(\l)\dagger}) \wedge d( \xi^{(\l)} \xi^{(\l)\dagger} ) \right\} \xi^{(\l)}  
= {} & d \xi^{(\l) \dagger} \wedge d \xi^{(\l)} + \xi^{(\l)\dagger} d \xi^{(\l)} \wedge \xi^{(\l)\dagger} d \xi^{(\l)} 
\nonumber \\ 
= {} & d A^{(\l)} + A^{(\l)} \wedge A^{(\l)} = F\ .
 \end{align}
On the other hand, since   $\big(\xi^{(\l)}_x\big| \h^{(\nu)}_x\big) \in \U_{m+k}$ for each $x \in \cU$,  we have 
\begin{multline}   
 \xi^{(\l) \dagger} \h^{(\nu)} = \xi^{(\l)\dagger}\left(\begin{array}{c} \nu \\ \Iden_{2k} \end{array} \right) (\Iden_{2k} + \nu^\dagger \nu)^{-\frac{1}{2} } = 0 \\
  \Longrightarrow\quad  \xi^{(\l)\dagger}\left(\begin{array}{c} \nu \\ \Iden_k \end{array} \right) = 0 \quad
 \text{and}\quad \Iden_{k + m} = \xi^{(\l)} \xi^{(\l)\dagger} + \h^{(\nu)} \h^{(\nu)\dagger} \ .
\end{multline} 
This implies that 
\begin{align*}   &   \xi^{(\l)\dagger} \left\{d (\xi^{(\l)} \xi^{(\l)\dagger}) \wedge ( \xi^{(\l)} \xi^{(\l)\dagger} ) \right\} \xi^{(\l)}  \\
&\qquad =  \xi^{(\l)\dagger} \left\{d \left(\Iden_{m+k} - \h^{(\nu)} \h^{(\nu)\dagger}\right) \wedge d \left(\Iden_{m+k} - \h^{(\nu)} \h^{(\nu)\dagger}\right) \right\}  \xi^{(\l)}    \\
&\qquad=   \xi^{(\l)\dagger} \left(d \left(\h^{(\nu)}  \h^{\nu) \dagger}\right) \wedge d \left(\h^{(\nu)}  \h^{(\nu)\dagger}\right) \right)  \xi^{(\l)} \\
&\qquad=    \xi^{(\l)\dagger} \left( \left(\begin{array}{c} d \nu \\ 0 \end{array} \right)\wedge (\Iden_{2k} + \nu^\dagger \nu)^{-1 }   \left(\begin{array}{cc} d \nu^\dagger & 0 \end{array} \right)\right)  \xi^{(\l)}   \ .
\end{align*}
From this and  \eqref{615} the conclusion follows.
\end{pf}
\par
\subsection{Prepotentials in canonical form}
We now  apply the above results to  instantons on a qk manifold $(M, g, Q)$ and its harmonic space $\HMqk$. 
Let    $\cW\subset  \HMqk$  be an open subset  of the harmonic space admitting   holomorphically adapted 
coordinates $(z^{\pm a}, u^i_{\pm})$ and   $G^o$ a compact Lie group  with complexification  $G$.   
By  compactness, we  may always assume that  $G^o \subset \U_m$ and $\gg^o = \Lie(G^0) \subset \gu_m$ for some  $m \geq 1$.  \par
We  denote  by   $(E^o, D^o)$  an instanton   on  $(M, g, Q)$   with structure group   $G^o$ and associated 
principal bundle with connection $(P^o, \o^o)$,  and  by  $(E'{}^o, D'{}^o)$ and $(P^o{}', \o'{}^o)$  the corresponding lifts  to   $\HMqk$.  
From  Theorem \ref{thm54},   Corollary \ref{lemma62}  and Lemma \ref{curvature1} we obtain
\begin{theo} 
Let $\l\colon \cW \simeq \cW {\times} \{I_2\}\to \bC_{k \times m} $ be a  smooth  map,   with  associated map $\xi^{(\l)}$ defined in \eqref{canonicalxi}  and $\nu {\,=} - \l^\dagger$,  such that for each $x \in  \cW {\times} \{I_2\}$ the following hold:
\beq \label{condlambda}
\begin{split}
\quad & \xi^{(\l)\dagger} \left(\frac{\p \xi^{(\l)} }{\p z^{+a}} + \frac{\p \xi^{(\l)}}{\p z^{-a}} \right)\bigg|_x \in \gg^o\ ,\\
& \frac{\p \l}{\p z^{\pm 1}} = \frac{\p \l}{\p z^{\pm 2}} = H_0{\cdot} \l = H_{++}{\cdot} \l = H_{--}{\cdot} \l = 0\ ,\\
&\left(d \nu\wedge   (\Iden_{2k} + \nu^\dagger \nu)^{-1 } d \nu^\dagger\right) \left(\frac{\p}{\p z^{\pm a}},  \frac{\p}{\p z^{\pm b}} \right)  =   0\ ,\\
&\left(d \nu\wedge   (\Iden_{2k} + \nu^\dagger \nu)^{-1 } d \nu^\dagger\right) \left(\frac{\p}{\p z^{+ a}},  \frac{\p}{\p z^{- b}} \right)  = \\
& \hskip 4 cm = - \left(d \nu\wedge   (\Iden_{2k} + \nu^\dagger \nu)^{-1 } d \nu^\dagger\right) \left(\frac{\p}{\p z^{- a}},  \frac{\p}{\p z^{+ b}} \right) .
\end{split}
\eeq
Then:
\begin{itemize}[leftmargin= 20pt]
\item[(1)]  For each $y \in \cW$, there exists a neighborhood $\cW' \subset \cW$ of $y$ on which there is   a  solution   $\mu^{(\l)}\colon \cW' \subset \cW\to \gg$  to the system of partial differential equations 
\beq \label{6.14-1} H_0{\cdot} \mu = 0\ ,\qquad \frac{\p \mu}{\p z^{-a}} = \xi^{(\l)\dagger} \frac{\p \xi^{(\l)}}{\p z^{-a}}\ , \qquad \frac{\p \mu}{\p z^{+1}} = \frac{\p \mu}{\p z^{+2}} = 0\,.\eeq
\item[(2)]  The    map 
\beq \label{U-prepotential1} 
A_{--}^{(\l)}\colon \cW' \to\gg\ ,\quad   A_{--}^{(\l)} \=  \Ad_{e^{-\mu}} \left( \xi^{(\l)\dagger}\big(H_{--}{\cdot}\xi^{(\l)}\big) + H_{--}{\cdot} \mu\right)  
\eeq
is    the prepotential   of an  instanton $(E^o|_{\cU}, D^o)$ on $\cU\=  \pi^{\operatorname{qk}}(\cW') \subset M$  with structure group $G^o$  in some analytic gauge.
\end{itemize}
Further,  any prepotential of an  instanton on $(M, g, Q)$  with structure group $G^o$ is gauge equivalent to  one of the form \eqref{U-prepotential1} for some $\l$  satisfying \eqref{condlambda} and  a solution $\mu$ of \eqref{6.14-1}. Finally, 
two  prepotentials   $A_{--}^{(\l)}$,   $A_{--}^{(\l')}$,    determined by   maps $\l, \l'\colon \cW' \to \bC_{k \times m}$,  are gauge   equivalent  
  if and only if the corresponding $1$-forms   $A^{o(\l)} =  \xi^{(\l)\dagger} d \xi^{(\l)}$  and $A^{o(\l')} =  \xi^{(\l')\dagger} d \xi^{(\l')}$ are equal (\footnote{This holds if and only if 
    $\l$ coincides with $\l'$  up to the  action of a row-permuting matrix in $\GL_{m+k}(\bC)$ as  mentioned in Corollary  \ref{lemma62}}).
\end{theo}
\begin{pf} Assume that $\l$ satisfies \eqref{condlambda} and set 
$$A^{o(\l)}\colon T^* \cW \to \bC_{m\times m}\ ,\qquad A^{o(\l)}|_x \= \xi^{(\l)\dagger} d \xi^{(\l)}\ .$$
By the first condition in \eqref{condlambda}, the map  $A^{o(\l)}$ is $\gg^o$ valued and  can be considered as the 
potential  of a connection $1$-form $\o^o$ on the trivial bundle $\pi'\colon P'{}^o = \cW \times G^o \to \cW$ and of 
the corresponding gauge field $(E', D')$.   By Lemma \ref{curvature1}  and Proposition \ref{prop28},  the other 
conditions in \eqref{condlambda} imply that   $(E', D')$   is the lift of an instanton on $\cU = \pi(\cW) \subset M$.  
 By Theorem \ref{thm53} for each $y \in \cW$  there  exists a neighborhood $\cW'$ of $y$  and an
analytic gauge for $P'|_{\cW'}= (P'{}^o|_{\cW'})^\bC$,  in which  the potential satisfies \eqref{claim1}.  This is tantamount to saying 
that there exists a solution to \eqref{6.14-1} on $\cW'$, i.e. that (1) holds. 
Claim (2) follows immediately from the fact that 
$A^{(\l)} = \Ad_{e^{-\mu}} \left( \xi^{(\l)\dagger}\big(H_{--}{\cdot}\xi^{(\l)}\big) + d \mu\right) $ is the potential of 
the lift of an instanton in an analytic gauge. \par
\smallskip
The property that any instanton $(E^o, D^o)$ admits a prepotential of the form \eqref{U-prepotential1} is checked as 
follows. By Corollary \ref{lemma62}, there exists a  gauge for the lift $(E'{}^o, D'{}^o)$ on $\HMqk$ in which the potential 
has the form  $\wt A^{(\l)} = \xi^{(\l)}{}^\dagger d\xi^{(\l)}$ for an appropriate choice of a matrix valued function $\l$.  
The map $\xi^{(\l)}$ which gives the potential in this gauge is  determined by pulling back the map 
$\cM\colon \cV_{m+k, k} \to \bC_{(m+k) \times k}$  in \S \ref{Stiefel1}  using  an appropriate $G^o$ equivariant 
embedding $\Phi'\colon P'{}^o \to \cV_{m+k, m}$ into the Stiefel bundle $\cV_{m+,k}$. 
Since $(E'{}^o, D'{}^o)$ is the lift of a gauge field on $M$, there is no loss of generality in  assuming that this equivariant 
embedding $\Phi'$  has the form $\Phi' = \Phi \circ \pi^{\text{qk}}$ for some $G_o$-equivariant embedding of 
$P^o$ into $\cV_{m+k, m}$. Under this hypothesis, the proof of Corollary \ref{lemma62} implies that the function $\l$ is 
constant along the fibres of the bundle $\pi^{\text{qk}}\colon \HMqk \to M$ and hence satisfies the first two lines of 
\eqref{condlambda}. The third line is immediate from  \eqref{curvature11} and the fact that  $\wt A^{(\l)}$  satisfies \eqref{3.10}. 
By Theorem \ref{thm53}, we know that there exists a gauge transformation $g = e^{\mu}\colon \cU \to G$ which modifies 
the above gauge of $P^o|_{\cW}$ (and  of its complexification $P|_{\cW}$) into  an analytic gauge. The new potential  given by such a gauge transformation has the component  $A_{--} = A(H_{--})$  as  in \eqref{U-prepotential1}.\par
\smallskip
For the final claim, we  observe that $A_{--}^{(\l)}$ and $A_{--}^{(\l')}$  are  respectively components of potentials 
which are  obtained from the  potentials $A^{o(\l)}$  and $A^{o(\l')}$ of the $G^o$-bundle  $P^o$ by  a gauge 
transformation in the complexification $P = (P^o)^\bC$.   Hence   $A_{--}^{(\l)}$ is gauge equivalent to  $A_{--}^{(\l')}$  
if and  only if $A^{o(\l)}$ is gauge equivalent to $A^{o(\l')}$.  By  Corollary \ref{lemma62}, this holds if and only if  
$A^{o(\l)} = A^{o(\l')}$. 
\end{pf}
We say that a prepotential of the form $A^{(\l)}_{--}$ is  {\it in canonical form} and  the maps  $\l$  and 
$\nu = - \l^\dagger$ are  its  {\it core} and  {\it conjugate core}, respectively.\par
\subsection{Examples of cores for prepotentials in canonical form on $\bH P^n$}
We now consider  the  basic example of qk manifold  $\bH P^n = \Sp_{n+1}/\Sp_n$, equipped with the  $\Sp_{n+1}$-invariant 
Riemannian metric  described in \S \ref{Hpn}. In this concluding section, we present  some   examples of cores $\l$, which determine prepotentials   for instantons on (open subsets of)  $\bH P^n$, with structure group  $G_o = \Sp_m$.  They correspond to the instantons on $\bH^n$ constructed in \cite{MS} and  are natural generalisations of some of the instantons  on $S^4 = \bH P^1$ appearing in the ADHM  classification \cite{At, ADHM}.  \par
\smallskip
From the discussion   in \S \ref{Hpn}, the Swann bundle and the harmonic space of $\bH P^n$  are 
$\bH^{n-1} {\setminus}\{0\}/\bZ_2$ and   $\left((\bH^{n+1} {\setminus}\{0\})/\bZ_2\right) \times \SL_2(\bC)$, respectively. 
Let us identify    $\cU_0 = \{ (q_0, \ldots, q_n) \in \bH^{n+1}\ :\ q_0 \neq 0\ \}  \subset \bH^{n+1}$ 
with  a subset of  $ \bC^{2n+2}$ using  the composition of the following two maps: First,
\beq
\label{6.14}  (q^0, q^1, \ldots, q^n)  \in \cU_0 \mapsto (q^0, r^1 \= q^1 (q^0)^{-1}, \ldots,r^n \=  q^n (q^0)^{-1}) \in  \bH^*\times  \bH^n ,
\eeq
 which is simply a local trivialisation of the  $\bH^*$-bundle $\pi\colon \bH^{n+1} {\setminus}\{0\} \to \bH P^n$
over $\cU_0 \subset \bH P^n$, and second 
\begin{multline}
\label{6.15}    \bH^* \times \bH^{n} \ni \left(q^0= \z^{0} + \z^{0' }{\bf j}; r^1 =  \z^{1} + \z^{1' }{\bf j} , \ldots,  r^n =  \z^{n} + \z^{n' }{\bf j} \right)  \\
 \longmapsto \left(\z^{0}, \z^{0' }, \z^{1} , \z^{1'}, \ldots ,  \z^{n},   \z^{n'} \right) \in  \bC^{2n+2} \ .
\end{multline}
This identification allows any point of  $\cU_0$
to be considered   as  a point of the  complexification $(\bC^{2n+2})^\bC \simeq \bC^2 \otimes \bC^{2n +2}$  using the embedding
\begin{multline*}
 \left(\z^0, \z^{0'}, \z^{1} ,  \z^{1'}, \ldots  \right)
 \\
\longmapsto \left(\begin{array}{llllllll} z^{1 1} \= \z^{0} &  z^{1 2} \= \z^{0'}  &  z^{1 3} \= \z^{1}  & z^{1 4} \= \z^{1'}  &  \ldots \\   
z^{2 1} \= - \overline{\z^{0'}}&  z^{2 2} \= \overline{\z^{0}}  &  z^{2 3} \= - \overline{\z^{1'}}  & z^{2 4} \= \overline{\z^{1}}  &   \ldots \\ 
\end{array}\right)  .
\end{multline*}
The coordinates on   $\cW = (\cU_0/\bZ_2) \times \SL_2(\bC) \subset \cH(\bH P^n)^{\text{qk}}$ which are  defined by 
\beq 
\left(x \equiv (z^{ia}), U= \left( u^i_\pm \right)\right)   \longmapsto \left(z^{\pm a} = u_i^\pm z^{i a}, (u_{\pm}^i) \right)
\ ,\quad u^\pm_i = \pm \e_{ij} u^j_\mp\ ,
\eeq
are   holomorphic adapted coordinates (as defined in \S \ref{sectionadapted}). Note that the   coordinates $(z^{i a}, u^j_{\pm})$ are  given in terms of   the  coordinates $(z^{\pm a}, u^i_{\pm})$ by   $z^{i a} = u^i_+ z^{+a} + u^i_{-} z^{-a}$ and that  $H_0{\cdot} z^{i a} = H_{++}{\cdot} z^{i a} = H_{--}{\cdot} z^{i a} = 0$  for each $i = 1,2$ and  $1 \leq a \leq 2n$. \par
\smallskip
For any quaternion $q = \z^{1} + \z^{2} {\bf j} \equiv (\z^1, \z^2)$, we  denote by $\bM(q)$ the   matrix representing the linear map   $R_q\colon \bH = \bC^2 \to \bH = \bC^2$ determined by the right action $R_q$  of $q$.  It is given by 
the $2 \times 2$  complex matrix
$$\bM(q) = \left(\begin{array}{cc} \z^{1} & - \overline{\z^{2}}\\ \z^{2} & \overline{\z^{1}}\end{array} \right) \ .$$
Note that the standard complex conjugate  $\overline{\bM(q)}$ of the matrix associated with a quaternion $q = x_0 + x_1 {\bf i} + x_2 {\bf j} + x_3 {\bf k}$,  is not equal to the 
matrix $\bM(\overline{q}^{\bH})$, associated with  the conjugate quaternion  $\overline{q}^{\bH} \= x_0 - x_1 {\bf i} - x_2 {\bf j} - x_3 {\bf k}$.  In fact, the proper relation is given by the conjugate transpose operator. Indeed, 
\beq 
\begin{split} \overline{\bM(\z^1 + \z^2 {\bf j})} = &\bM(\overline{\z^1} + \overline{\z^2} {\bf j})\ ,\qquad \bM(\z^1 + \z^2 {\bf j})^T =  \bM(\z^1 -  \overline{\z^2} {\bf j})\ ,\\
& \bM(\z^1 + \z^2 {\bf j})^\dagger = \bM(\overline{\z^1} - \z^2 {\bf j}) = \bM(\overline{(\z^1 + \z^2 {\bf j})}^{\bH})\ .
\end{split}
\eeq
\par \smallskip
We now pick two  integers $k, m \geq 1$ and consider an $(n+1)$-tuple of  complex matrices $\underset{(a)} \cA$, $0 \leq a \leq n$,  in $\bC_{2m \times 2k}$,   having the  form  
\beq \label{1} \underset{(a)} \cA = \left(\begin{array}{cccc} \underset{(a)} \cA{}_{1}^1 & \underset{(a)} \cA{}^1_{2} & \ldots & \underset{(a)} \cA{}^1_{k}\\
\underset{(a)}\cA{}_{1}^2 & \underset{(a)}\cA{}^2_2 & \ldots & \underset{(a)}\cA{}^2_k\\
 \vdots & \vdots & \ddots &  \vdots\\
 \underset{(a)}\cA{}_{1}^k & \underset{(a)}\cA{}^k_2 & \ldots & \underset{(a)}\cA{}^k_k \end{array}\right)\ \ \text{with $2\times 2$ blocks}\ \ \underset{(a)} \cA{}^i_j = \bM(\underset{(a)}q{}^i_j), \  \underset{(a)}q{}^i_j \in \bH\ .\eeq
Note that by identifying each block  $\bM(\underset{(a)}q{}^i_j)$ with the corresponding quaternion $q^i_j$, the matrix $\underset{(a)} \cA $ can be considered as an element of $\bH_{m \times k}$
on which the   $({\cdot})^\dagger$  operator acts as   the  $\bH$-conjugate transpose operator.
\par
Finally, we  consider  the map $ \nu\colon \cU_0 \to  \bC_{2m \times 2 k}$ defined by 
  \begin{multline} \label{nuzed}  \mathfrak \nu|_{(z^{j a})} \= \\ 
  \underset{(0)} \cA   + \sum_{a = 1}^{n} \underset{(a)}\cA \left( \begin{array} {cccc} \smallmatrix \bM(z^{1(2a)} {+} z^{1(2 a{+}1)}{\bf j})  \endsmallmatrix& 0 & \ldots & 0\\
  0 & \smallmatrix\bM(z^{1(2a)} {+} z^{1  (2a{+}1)}{\bf j})  \endsmallmatrix& \ldots & 0\\
  \vdots & \vdots  & \ddots & \vdots\\
  0 & 0  & \ldots & \smallmatrix\bM(z^{1 (2a)} {+}  z^{1  (2a{+}1)}{\bf j})  \endsmallmatrix
  \end{array}\right)   \end{multline}
and we set 
$\l = - \nu^\dagger$.
For this map   $\l$ we have:
\begin{itemize} [itemsep=4pt plus 4pt minus 2pt, leftmargin=18pt,topsep=0pt]
\item  At each point the corresponding map   $\xi^{(\l)}$ is  equal to a matrix with $2\times 2$ blocks of the form $\bM(q)$ for some quaternion $q$. It is therefore identifiable with a quaternionic matrix in $\bH_{(k + m) \times m}$. This implies that for each $x \in \cU_0 \simeq \cW \times \{I_2\}$, the matrix 
 $\xi^{(\l)\dagger} \left(\frac{\p \xi^{(\l)} }{\p z^{+a}} + \frac{\p \xi^{(\l)}}{\p z^{-a}} \right)\Big|_x$
  is an element in $\su_{2 m} {\cap} \ggl_m(\bH) = \sp_m$, i.e.  $\l$ satisfies  $\eqref{condlambda}_{I}$ with $\gg^0 = \sp_m$
\item The  matrix valued function $\l_{(z^{ia})}$ is independent of the coordinates $z^{\pm 1}$, $z^{\pm 2}$ and has trivial directional derivative  along each of the vector fields $H_A$, $A \in \{0, ++, --\}$. This implies that $\l$ trivially satisfies  $\eqref{condlambda}_{II}$.
\end{itemize}
Due to this, if $\nu$ is such that  conditions $\eqref{condlambda}_{III}$ e $\eqref{condlambda}_{IV}$ are   also satisfied,
 we may conclude that $\l = - \nu^\dagger$ is the core of a prepotential. 
We now claim   that   these two conditions are  satisfied  if we  impose  that  the  $\underset{(a)} \cA$  satisfy the additional  algebraic conditions
\beq  \label{nuzed2}  \underset{(a)}\cA   \underset{(a)}\cA^\dagger - \overline{\underset{(a)}\cA   \underset{(a)}\cA^\dagger}^\bH = 0\ ,\qquad  \underset{(a)}\cA  \underset{(b)}\cA^\dagger =  \underset{(b)}\cA  \underset{(a)}\cA^\dagger \ ,\qquad 0 \leq a, b \leq n\ .\eeq
In fact, these conditions  imply  that  the quaternionic entries of $\Iden_{2k} + \nu^\dagger \nu$  
 have  zero imaginary parts, that is they are real.  Thus, for any $\a, \b \in \{+, -\}$ and $2 \leq a,b \leq  n+1$, we have 
\begin{multline*} \left(d \nu\wedge   (\Iden_{2k} + \nu^\dagger \nu)^{-1 } d \nu^\dagger\right)\left(\frac{\p}{\p z^{\a a}},  \frac{\p}{\p z^{\b b}} \right)  = \\
= \left\{\begin{array}{cc} {\pm} \left( \underset{(a)}\cA  (\Iden_{2k} + \nu^\dagger \nu)^{-1 } \underset{(b)}\cA^\dagger\right) \!\!(u_\a^1 u_\b^2 {-} u_\b^1 u_\a^2) & \text{if}\ a, b\ \text{are both even or odd}\\
\left( \underset{(a)}\cA  (\Iden_{2k} + \nu^\dagger \nu)^{-1 }\bM({\bf j}) \underset{(b)}  \cA^\dagger \right)\!\!(u_\a^1 u_\b^2 {-} u_\b^1 u_\a^2) & \text{if}\ a\ \text{is even and}\ b\ \text{is odd}.
 \end{array}
\right.  
\end{multline*}
 From this  $\eqref{condlambda}_{III}$ and $\eqref{condlambda}_{IV}$ follow, as  claimed. 
\par
\medskip
We  conclude that for any choice of the matrices as in \eqref{1},  \eqref{nuzed} 
and satisfying  \eqref{nuzed2},  the corresponding map $\l = - \nu^\dagger$ is the core of a prepotential of an instanton 
on $\wh \cU_0 = \pi(\cU_0) \subset \bH P^n$ with structure group $G^o = \Sp_m$.  As mentioned above, instantons 
constructed in this way are  precisely  those  constructed   in \cite[\S 3]{MS}, where it was also proved that  they are  
well defined  over the entire  $\bH P^n$ and are  characterised by certain   topological properties.  Further details of 
this class of instantons in the lowest dimensional case of  $ \bH P^1 = S^4$ can be found for instance in  \cite[Ch. II]{At}.   
Explicit coordinate expressions for    prepotentials for instantons on $S^4$ with structure group $G^o = \Sp_1$ were determined   in  \cite{Og}.
\par
\medskip

\vskip 1.5truecm
\hbox{\parindent=0pt\parskip=0pt
\vbox{\baselineskip 9.5 pt \hsize=3.1truein
\obeylines
{\smallsmc
Chandrashekar Devchand
Max-Planck-Institut f\"ur Gravitationsphysik 
(Albert-Einstein-Institut)
Am M\"uhlenberg 1 
D-14476 Potsdam 
Germany
}\smallskip
{\smallit E-mail}\/: {\smalltt devchand@aei.mpg.de
}
}
\hskip  -3truemm
\vbox{\baselineskip 9.5 pt \hsize=3.7truein
\obeylines
{\smallsmc
Massimiliano Pontecorvo
Dipartimento di Matematica e Fisica
Universit\`a Roma Tre
Largo San Leonardo Murialdo,1
I-00146 Roma
Italy
}\smallskip
{\smallit E-mail}\/: {\smalltt max@mat.uniroma3.it}
}
}
\vskip 1.0 truecm
\hbox{\parindent=0pt\parskip=0pt
\vbox{\baselineskip 9.5 pt \hsize=3.1truein
\obeylines
{\smallsmc
Andrea Spiro
Scuola di Scienze e Tecnologie
Universit\`a di Camerino
Via Madonna delle Carceri 
I-62032 Camerino (Macerata)
Italy
}\smallskip
{\smallit E-mail}\/: {\smalltt andrea.spiro@unicam.it
}
}
}
\end{document}